\documentclass[a4paper,11pt,twoside,reqno]{amsart}

\theoremstyle{plain}
\newtheorem{theorem}                {Theorem}      [section]
\newtheorem*{theorem*}                {Theorem \ref{thm:appl}}
\newtheorem{proposition}  [theorem]  {Proposition}

\newtheorem{lemma}        [theorem]  {Lemma}

\theoremstyle{definition}

\newtheorem{remark}       [theorem]  {Remark}

\numberwithin{equation}{section}

\DeclareMathOperator{\trace}{trace}
\DeclareMathOperator{\grad}{grad}
\DeclareMathOperator{\Span}{span}

\DeclareMathOperator{\sgn}{sgn}
\DeclareMathOperator{\Sol}{Sol}
\DeclareMathOperator{\Id}{Id}

\usepackage[colorlinks=true,linkcolor=blue,citecolor=blue]{hyperref}

\begin{document}
	
	\title[Biconservative surfaces in the $4$-dimensional hyperbolic space]
	{Extrinsic characterizations of biconservative surfaces in the $4$-dimensional hyperbolic space}
	
	\author{Simona Nistor, Mihaela Rusu}
	
	\address{Faculty of Mathematics, Al. I. Cuza University of Iasi,
		Blvd. Carol I, 11 \\ 700506 Iasi, Romania} \email{nistor.simona@ymail.com}
	
	\address{Faculty of Mathematics, Al. I. Cuza University of Iasi,
		Blvd. Carol I, 11 \\ 700506 Iasi, Romania} \email{mihaelarussu10@yahoo.com}
	
	\begin{abstract}
		Biconservative submanifolds arise as a natural relaxation of the biharmonic condition and play an important role in the submanifold theory. In this paper, we study non-CMC biconservative surfaces with parallel normalized mean curvature vector field (PNMC surfaces) in the four-dimensional hyperbolic space $\mathbb{H}^4$, for which we consider the hyperboloid model. We provide a local extrinsic description of such surfaces, showing that they are generated by a directrix curve lying in a totally geodesic hypersurface $\mathbb{H}^3$ of $\mathbb{H}^4$, through a certain normal flow. This extrinsic classification of non-CMC, PNMC biconservative surfaces in $\mathbb{H}^4$ splits naturally into three cases according to the type of a certain vector field, which can be non-zero null, spacelike or timelike. We also prove that these surfaces are invariant under the action of a parabolic, elliptic, and hyperbolic one-parameter group of isometries of $\mathbb{H}^4$, respectively. Moreover, their full groups of ambient isometries preserving the surfaces are determined. Together with the previous results, the classification of non-CMC, PNMC surfaces in four-dimensional space forms is now complete, from both intrinsic and extrinsic points of view.
	\end{abstract}
	
	
	
	\maketitle
	
	\section{Introduction}
	Biconservative submanifolds have become an important topic in the modern study of submanifold geometry, as they can be viewed as a natural relaxation of the biharmonic condition.
	%Instead of requiring the bitension field to vanish, one asks only that its tangential component to be zero (see, for example, \cite{AK, CMOP, Jiang87, NPhD17}).
	Biharmonic submanifolds generalize the classical notion of minimal submanifolds and are defined as isometric immersions $\varphi:\left(M^m,g\right)\to\left(N^n,h\right)$ satisfying the biharmonic equation
	$$
	\tau_2\left(\varphi\right)=-\Delta^\varphi \tau\left(\varphi\right)-\trace R^N\left(\varphi_\ast,\tau\left(\varphi\right)\right)\varphi_\ast=0.
	$$
	Here, $\Delta^\varphi$ denotes the rough Laplacian acting on sections of the pull-back bundle $\varphi^{-1}\left(TN\right)$, $R^N$ is the curvature tensor field of $N$, and
	$$
	\tau\left(\varphi\right)=mH
	$$
	is the tension field associated with $\varphi$, where $H$ is the mean curvature vector field. For an overview of the current research on biharmonic submanifolds, we refer to \cite{AN, FO2022, OC-20}. The bitension field $\tau_2\left(\varphi\right)$ decomposes into tangential and normal components, and biconservative submanifolds are characterized by the vanishing of its tangential part. For further details on the geometric meaning of the equation $\tau_2\left(\varphi\right)^\top=0$, see, for example, \cite{AK, CMOP, Jiang87, LMO, NPhD17}.
	
	Numerous results have already been obtained related to biconservative submanifolds in various ambient spaces, including Euclidean spaces, Euclidean spheres, hyperbolic spaces, product spaces, $\Sol_3$, de Sitter spaces, among others (see, for example \cite{F2025, FOP2015, F2015, Kayhan2025, MOP2023, MOR2016JGA, Turgay2015, TurgayYegin2025}). 
	
	The study of biconservative submanifolds was initially focused on the case of hypersurfaces in space forms, i.e., in spaces with constant sectional curvature (see, for example, \cite{CMOP, Hasanis-Vlachos}). Then, it was extended to the investigation of biconservative submanifolds of codimension two in space forms, with particular emphasis on the case of surfaces. In this setting, surfaces with parallel mean curvature vector field (PMC surfaces) turn out to be trivially biconservative; therefore, the most interesting situation arises when considering non-PMC surfaces in four-dimensional space forms, namely in $N^4\left(\varepsilon\right)$. 
	
	Within this class of surfaces, one may distinguish two subclasses: surfaces with constant mean curvature vector field (CMC surfaces) and non-CMC surfaces. The classification of CMC biconservative surfaces in four-dimensional space forms was obtained in \cite{MOR2016JGA}. In the non-CMC case, additional geometric assumptions are required in order to derive classification results. A particularly useful condition, which considerably simplifies the structure equations and has played a crucial role in obtaining classification results, is for the surface to have parallel normalized mean curvature vector field in the normal bundle, i.e., to be a PNMC surface. We mention here that the substantial codimension of a non-CMC, PNMC biconservative surface in a space form $N^n\left(\varepsilon\right)$, $n\geq 5$, is two (see \cite{NOTYS, SNMR2024, YeginTurgay2018} and also \cite{YeginTurgay2025}). The complementary non-PNMC case has also been investigated under additional geometric assumptions (see \cite{AMOS}).
	
	Non-CMC, PNMC biconservative surfaces in $N^4\left(\varepsilon\right)$, for $\varepsilon=0$ and $\varepsilon=1$, i.e., in Euclidean spaces and Euclidean spheres, respectively, were studied in \cite{YeginTurgay2018} and \cite{NOTYS} from both intrinsic and extrinsic points of view. In hyperbolic ambient spaces, i.e., when $\varepsilon=-1$, obtaining extrinsic descriptions of biconservative surfaces remains a delicate problem. In $\mathbb{H}^4$, the biconservative condition interacts subtly with the Lorentzian structure of the ambient space $\mathbb{R}^5_1$ in the hyperboloid model, and a full extrinsic characterization is far from immediate.
	
	The present work continues the intrinsic investigation initiated in \cite{SNMR2024}, where non-CMC, PNMC biconservative surfaces in $\mathbb{H}^4$ were studied from an intrinsic point of view. In particular, it was proved that they form a two-parameter family. Here, we complement that analysis by providing a full extrinsic description of such surfaces based on the geometry of the ambient space $\mathbb{R}^5_1$. 
	
	As a preliminary step, we establish two auxiliary results for curves in $\mathbb{H}^3$, relating their curvature and torsion to the existence of a constant vector satisfying a prescribed inner-product condition with the position vector of the curve (Lemma \ref{lemma-1} and Lemma \ref{lemma-2}). These results will be used to establish the compatibility of the additional inner-product condition with the prescribed curvature and torsion of the directrix curve appearing in the local parametrizations.
	
	We then show that any non-CMC, PNMC biconservative surface in $\mathbb{H}^4$ admits a local description in terms of a directrix curve $\sigma$, which lies in a totally geodesic hypersurface $\mathbb{H}^3$ of $\mathbb{H}^4$, and of generating curves given by the integral curves $\gamma$ of a prescribed vector field. The directrix curve is determined by the mean curvature function $f$, which satisfies a specific second-order ordinary differential equation. Our classification naturally splits into three distinct cases, depending on the sign of an integration constant $C$, which gives rise to an associated vector field \(\xi\) whose causal type is determined by the sign of \(C\). In the first case, when $C=0$, the vector field $\xi$ is null and the generating curves $\gamma$ of the surface are parabolas contained in lightlike affine planes (Theorem \ref{theorem-C=0}). In the second case, when $C>0$, the vector field $\xi$ is spacelike and the generating curves $\gamma$ are circles contained in Euclidean affine planes (Theorem \ref{theorem-C>0}). In the last case, when $C<0$, the vector field $\xi$ is timelike and the generating curves $\gamma$ are branches of hyperbolas contained in Lorentzian affine planes (Theorem \ref{theorem-C<0}). This distinction also leads to three different geometric descriptions of the directrix curve in $\mathbb{H}^3$: in the three cases, the directrix meets transversely a foliation of $\mathbb{H}^3$ by horospheres, equidistant surfaces, and geodesic spheres, respectively. 
	
	We also determine the isometries of $\mathbb{H}^4$ preserving these surfaces. According to the sign of $C$, the surfaces admit a one-parameter group of parabolic, elliptic, or hyperbolic rotational symmetries for $C=0$, $C>0$, or $C<0$, respectively. In each case, the full group of isometries of $\mathbb{H}^4$ preserving the surface is generated by this one-parameter group together with a Lorentzian reflection.
	
	Together with the analogous results in the Euclidean and spherical settings, this paper completes the local classification of non-CMC, PNMC biconservative surfaces in four-dimensional space forms and provides a foundation for further investigations of their global properties. Our approach also illustrates the interplay between intrinsic geometry and extrinsic geometry, a recurring theme in the broader study of submanifold theory.
	
	\textbf{Conventions and notations.}
	In general, all Riemannian metrics are denoted by the same symbol $\langle\cdot,\cdot\rangle$. When no confusion arises, we omit explicit reference to the metric. All manifolds are assumed to be connected and oriented. For the rough Laplacian acting on sections of the pull-back bundle $\varphi^{-1}\left(TN^n\right)$ and for the curvature tensor field, we adopt the following sign conventions:
	$$
	\Delta^{\varphi}=-\trace\left(\nabla^{\varphi}\nabla^{\varphi}-\nabla^{\varphi}_{\nabla}\right)
	$$
	and
	$$
	R\left(X,Y\right)Z=[\nabla_X,\nabla_Y]Z-\nabla_{[X,Y]}Z,
	$$
	respectively. Here, $\varphi:M^m\to N^n$ is a smooth map between two Riemannian manifolds, $\nabla^\varphi$ denotes the induced connection on $\varphi^{-1}\left(TN^n\right)$, and $\nabla$ is the Levi-Civita connection on $M$. Moreover, we denote by $\hat{\nabla}$ and $\tilde{\nabla}$ the Levi-Civita connections of the Lorentz-Minkowski space $\mathbb{R}^{m+1}_1$ and of the hyperbolic space $\mathbb{H}^m$, endowed with the standard metrics, respectively.
	
	Our approach is essentially local. In order to avoid trivial cases in the study of non-CMC, PNMC biconservative surfaces $M^2$ in $\mathbb{H}^4$, we assume that the mean curvature function of the surface is positive, its gradient is nowhere vanishing, $\nabla^\perp \left({H/|H|}\right)=0$, and $M$ is completely contained in $\mathbb{H}^4$, in the sense that any open subset of $M$ cannot lie in a totally geodesic hypersurface $\mathbb{H}^3\subset\mathbb{H}^4$. For simplicity, throughout the paper, whenever we refer to a PNMC biconservative immersion or surface in $\mathbb{H}^4$, all these assumptions are understood to hold.
	
	\textbf{Acknowledgments.} The authors would like to thank Cezar Oniciuc for many valuable discussions and comments.
	
	\section{An extrinsic approach to PNMC biconservative surfaces in $\mathbb{H}^4$}
	%%%%%%%%%%%%%%%%%%%%%%%%%%%%%%%%%%%%%%%%%%%%%%%%%%%%%%%%%%%%%%%%%%%%%%
	For the four-dimensional hyperbolic space $\mathbb{H}^4$ we consider the hyperboloid model. More precisely, in the Lorentz-Minkowski space $\mathbb{R}_1^5=\left(\mathbb{R}^5,\langle \cdot , \cdot \rangle\right)$, we denote by $\langle \cdot , \cdot\rangle$  the bilinear form
	$$
	\langle \overline{x}, \overline{y}\rangle=\sum_{i=1}^4 x^i y^i-x^5 y^5,
	$$
	and thus,
	$$
	\mathbb{H}^4=\left\{\overline{x} \in \mathbb{R}_1^5:\langle\bar{x}, \bar{x}\rangle=-1 \text { and } x^5>0\right\} \text {, }
	$$
	i.e., $\mathbb{H}^4$ is the upper part of the hyperboloid of two sheets.
	
	We first present two results concerning geometric properties of curves in $\mathbb{H}^3$, when $\mathbb{H}^3$ may be seen as a totally geodesic hypersurface of $\mathbb{H}^4$.
	\begin{lemma}\label{lemma-1}
		Let $\sigma:I\to\mathbb{H}^3$ be a curve parametrized by arc-length, with curvature $\kappa>0$ and torsion $\tau>0$, and let $\{V_1,V_2,V_3\}$ be its Frenet frame. Assume that there exists a constant vector $\overline{b}_2\in\mathbb{R}^4_1$ such that
		\begin{equation}\label{prop-b2-lemma1}
			\langle \sigma\left(u\right),\overline{b}_2\rangle=\Theta,
			\qquad
			\langle\overline{b}_2,\overline{b}_2\rangle=\delta,
			\qquad\text{and}\qquad
			\langle\overline{b}_2,V_3\rangle\neq 0,
		\end{equation}
		where $\Theta\in C^\infty\left(I\right)$ and $\delta\in\{-1,0,1\}$. Then
		\begin{equation}\label{positive-nominator}
			\delta+\Theta^2-\left(\Theta'\right)^2
			-\left(\frac{\Theta''-\Theta}{\kappa}\right)^2>0,
		\end{equation}
		and the torsion of $\sigma$ is uniquely determined by $\kappa$, $\Theta$, and $\delta$ by
		\begin{equation}\label{general-tau}
			\tau=
			\frac{
				\left|
				\Theta'''+\left(\kappa^2-1\right)\Theta'
				-\frac{\kappa'}{\kappa}\left(\Theta''-\Theta\right)
				\right|
			}{
				\kappa
				\sqrt{
					\delta+\Theta^2-\left(\Theta'\right)^2
					-\left(\frac{\Theta''-\Theta}{\kappa}\right)^2
				}
			}.
		\end{equation}
	\end{lemma}
	\begin{proof}
		We express the constant vector $\overline{b}_2\in\mathbb{R}^4_1$ with respect to the pseudo-orthonormal frame $\{\sigma,V_1,V_2,V_3\}$ along $\sigma$ as
		\begin{equation}\label{b2-coef}
			\overline{b}_2=a_0\sigma+a_1V_1+a_2V_2+a_3V_3,
		\end{equation}
		where $a_i\in C^\infty\left(I\right)$, $i\in\{0,\ldots,3\}$.
		
		Since $\sigma$ is a curve in $\mathbb{H}^3$, we have
		$\langle\sigma,\sigma\rangle=-1$. Taking the inner product of
		\eqref{b2-coef} with $\sigma$ and using
		$\langle\sigma,\overline{b}_2\rangle=\Theta$, we obtain
		\begin{equation}\label{expression-a0}
			a_0=-\Theta.
		\end{equation}
		Differentiating
		$\langle\sigma,\overline{b}_2\rangle=\Theta$ and using $\sigma'=V_1$ and the fact that $\overline{b}_2$ is a constant vector, we obtain
		\begin{equation}\label{expression-a1}
			a_1=\langle\overline{b}_2,V_1\rangle=\Theta'.
		\end{equation}
		Next, differentiating \eqref{expression-a1} and using the Frenet formula
		$$
		\hat{\nabla}_{V_1}V_1=\kappa V_2+\sigma,
		$$
		we obtain
		$$
		\Theta''
		=
		\langle\overline{b}_2,\hat{\nabla}_{V_1}V_1\rangle
		=
		\kappa\langle\overline{b}_2,V_2\rangle
		+
		\langle\overline{b}_2,\sigma\rangle.
		$$
		Therefore,
		\begin{equation}\label{expression-a2}
			a_2=\langle\overline{b}_2,V_2\rangle
			=\frac{\Theta''-\Theta}{\kappa}.
		\end{equation}
		It remains to determine $a_3$. Since
		$\langle\overline{b}_2,\overline{b}_2\rangle=\delta$ and
		$a_3=\langle\overline{b}_2,V_3\rangle\neq0$, from
		\eqref{b2-coef}--\eqref{expression-a2} we obtain
		\begin{equation}\label{expression-a3-squared}
			a_3^2=
			\delta+\Theta^2-\left(\Theta'\right)^2
			-\frac{\left(\Theta''-\Theta\right)^2}{\kappa^2}
			> 0.
		\end{equation}
		Thus, \eqref{positive-nominator} follows.
		
		Finally, differentiating \eqref{expression-a2} and using the Frenet formula
		$$
		\hat{\nabla}_{V_1}V_2=-\kappa V_1+\tau V_3,
		$$
		we obtain
		$$
		a_2'
		=
		-\kappa\langle\overline{b}_2,V_1\rangle
		+\tau\langle\overline{b}_2,V_3\rangle
		=
		-\kappa\Theta'+\tau a_3.
		$$
		Hence,
		\begin{equation}\label{expression-a3}
			a_3=\frac{a_2'+\kappa\Theta'}{\tau}
			=
			\frac{
				\Theta'''+\left(\kappa^2-1\right)\Theta'
				-\frac{\kappa'}{\kappa}\left(\Theta''-\Theta\right)
			}{
				\kappa\tau
			}
			\neq0.
		\end{equation}
		Combining \eqref{expression-a3-squared} and \eqref{expression-a3}, and using $\kappa>0$ and $\tau>0$, we obtain \eqref{general-tau}.
	\end{proof}
	A direct computation also yields the converse of the previous result. More precisely, we have the following.
	\begin{lemma}\label{lemma-2}
		Let $\sigma:I\to\mathbb{H}^3$ be a curve parametrized by arc-length, with curvature $\kappa>0$ and torsion $\tau>0$, and let $\{V_1,V_2,V_3\}$ be its Frenet frame. Let $\Theta\in C^\infty\left(I\right)$ and $\delta\in\{-1,0,1\}$. Assume that \eqref{positive-nominator} and \eqref{general-tau} hold. Then there exists a unique constant vector $\overline{b}_2\in\mathbb{R}^4_1$, given by
		\begin{equation}\label{expression-b2}
			\begin{aligned}
				\overline{b}_2=-\Theta\sigma+\Theta'V_1
				+\frac{\Theta''-\Theta}{\kappa}V_2\
				+\frac{
					\Theta'''+\left(\kappa^2-1\right)\Theta'
					-\frac{\kappa'}{\kappa}\left(\Theta''-\Theta\right)
				}{
					\kappa\tau
				}V_3,
			\end{aligned}
		\end{equation}
		which satisfies \eqref{prop-b2-lemma1}.
	\end{lemma}
	\begin{remark}\label{remark-b2-H3}
		Assume that $\delta=-1$ and that $\Theta<0$. Then
		$\langle\overline{b}_2,\overline{b}_2\rangle=-1$,
		so $\overline{b}_2$ belongs to one of the two sheets of the hyperboloid in $\mathbb{R}^4_1$. On the other hand, since $\sigma\left(u\right)\in\mathbb{H}^3$ and
		$$
		\langle\sigma\left(u\right),\overline{b}_2\rangle
		=
		\Theta\left(u\right)<0,
		$$
		the timelike vectors $\sigma\left(u\right)$ and $\overline{b}_2$ have the same time orientation. Hence, $\overline{b}_2$ belongs to the upper part of the hyperboloid as $\sigma\left(u\right)$, and therefore, $\overline{b}_2\in\mathbb{H}^3$.
	\end{remark}
	Now, let $\varphi: M^2 \to \mathbb{H}^4$ be a PNMC biconservative immersion, $\iota: \mathbb{H}^4 \to \mathbb{R}^5_1$ be the canonical inclusion, and denote
	$$
	\Phi=\iota \circ \varphi: M \rightarrow \mathbb{R}^5_1.
	$$
	Recall that $N M$ represents the normal bundle of $\varphi$ and, in order to avoid any confusion, we denote by $N_{\Phi} M$ the normal bundle of the immersion $\Phi$. Clearly, the two normal bundles are related by
	$$
	N_{\Phi} M=\iota_*\left(N M^2\right) \oplus \operatorname{span}\{\Phi\}
	$$
	and we have
	$$
	B_{\Phi}\left(X, Y\right)=\iota_*\left(B\left(X, Y\right)\right)+\langle X, Y \rangle \Phi, \quad X, Y \in C\left(T M\right)
	$$
	where $B_{\Phi}$ denotes the second fundamental form of $\Phi$.
	
	As usual, when we work with isometric immersions, we identify $M$ with $\varphi\left(M\right)$ or $\Phi\left(M\right)$ and a tangent vector field $X \in C\left(TM\right)$ with $\varphi_* X$ or $\Phi_* X$.
	
	The mean curvature vector field of $M$ in $\mathbb{H}^4$ is defined by $H=\left(\trace B\right)/2\in C\left(NM\right)$, where the $\trace$ is considered with respect to the domain metric. The mean curvature function is defined by $f=|H|$.
	
	According to our conventions, $f$ is a smooth positive function on $M$, and we can define the vector fields
	$$
	E_1=\frac{\grad f}{|\grad f|} \qquad \text{and} \qquad E_3=\frac{H}{f}.
	$$
	Moreover, if we denote by $A_{3}$ the shape operator of $M$ corresponding to $E_3$, it follows that $M$ is a PNMC biconservative surface in $\mathbb{H}^4$, i.e.,  $\tau_2^\top\left(\varphi\right)=0$, if and only if 
	\begin{equation}\label{eq:PNMC-biconservative}
		A_{3}\left(\grad f\right)=-f\grad f.
	\end{equation}
	(see \cite{FLO2021}). Because of orientation, we can consider the positively oriented global orthonormal frame fields $\left\{E_1,E_2\right\}$ in the tangent bundle $TM$ and $\left\{E_3,E_4\right\}$ in the normal bundle $NM$.
	
	Clearly, $E_2f=0$. Denoting by $A_4$ the shape operator corresponding to $E_4$, we recall the following result from \cite{SNMR2024}.
	\begin{theorem}[\cite{SNMR2024}]\label{thm:fundamentalProperties}
		Let $\varphi:M^2\to\mathbb{H}^4$ be a PNMC biconservative surface. Then, the following hold:
		\begin{itemize}
			\item [(i)] the Levi-Civita connection $\nabla$ of $M$ and the normal connection $\nabla^\perp$ of $M$ in $\mathbb{H}^4$ are given by
			\begin{equation}\label{Levi-Civita-connection-f}
				\nabla_{E_1}E_1=\nabla_{E_1}E_2=0, \quad \nabla_{E_2}E_1=-\frac{3}{4}\frac{E_1 f}{f}E_2, \quad \nabla_{E_2}E_2=\frac{3}{4}\frac{E_1 f}{f}E_1
			\end{equation}
			and
			\begin{equation*}\label{normal-connection-f}
				\nabla^\perp E_3=0, \qquad \nabla^\perp E_4=0;
			\end{equation*}
			\item [(ii)] the shape operators corresponding to $E_3$ and $E_4$ are given, with respect to $\left\{E_1, E_2\right\}$, by the matrices
			\begin{equation*}\label{shape-operators-A3-A4-f}
				A_3=\left(
				\begin{array}{cc}
					-f & 0 \\
					0 & 3f
				\end{array}
				\right),\quad
				A_4=\left(
				\begin{array}{cc}
					cf^{3/2} & 0 \\
					0 & -cf^{3/2}
				\end{array}
				\right),
			\end{equation*}
			where $c$ is some non-zero real constant;
			\item [(iii)] the Gaussian curvature $K$ and the mean curvature function $f$ are related by
			\begin{equation}\label{relation-K-f}
				K=-1-3f^2-c^2f^3,
			\end{equation}
			thus $1+K<0$ on $M$;
			\item [(iv)] the mean curvature function $f$ satisfies
			\begin{equation}\label{second-order-invariant-f}
				f\Delta f+\left|\grad f\right|^2-\frac{4}{3}f^2-4f^4-\frac{4}{3}c^2f^5=0.
			\end{equation}
		\end{itemize}
	\end{theorem}
	\begin{remark}
		The constant $c$ which appears in the expression of the shape operators is uniquely determined by the metric induced on $M$ by the immersion $\varphi$ and it is not an indexing parameter. In \cite{SNMR2024}, it was proved that if $\left(M^2,g\right)$ is an abstract surface, then it admits at most one PNMC biconservative immersion in $\mathbb{H}^4$, up to isometries of $\mathbb{H}^4$. Moreover, the set of all abstract surfaces $\left(M^2,g\right)$ that admit a (unique) PNMC biconservative immersion in $\mathbb{H}^4$ forms a two-parameter family, indexed by $c$ and $C$, where the real constant $C$ is given by the first integral of the second-order ordinary equation determined by \eqref{second-order-invariant-f}.
	\end{remark}
	From Theorem \ref{thm:fundamentalProperties}, by straightforward computations, we get 
	\begin{equation}\label{conditions-f}
		E_2\left(E_1f\right)=0\qquad\text{and}\qquad E_2\left(E_1\left(E_1f\right)\right)=0.
	\end{equation} 
	Then, we can see that the second fundamental form of $M$ in $\mathbb{H}^4$ and the Levi-Civita connection of $M$ in $\mathbb{R}^5_1$ are given by
	\begin{equation*}
		B\left(E_1,E_1\right)=-fE_3+cf^{3/2}E_4, \qquad	B\left(E_1,E_2\right)=0, \qquad B\left(E_2,E_2\right)=3fE_3-cf^{3/2}E_4,
	\end{equation*}
	and 
	\begin{equation}\label{formulas Th.3.1}
		\left\{
		\begin{aligned}
			&\hat{\nabla}_{E_1} E_1=\tilde{\nabla}_{E_1} E_1+\Phi=\nabla_{E_1}E_1+B\left(E_1,E_1\right)+\Phi=-f E_3+c f^{3 / 2} E_4+\Phi \\
			&\hat{\nabla}_{E_2} E_1=\tilde{\nabla}_{E_2} E_1=\nabla_{E_2}E_1+B\left(E_2,E_1\right)=-\frac{3}{4} \frac{E_1 f}{f} E_2 \\
			&\hat{\nabla}_{E_1} E_2=\tilde{\nabla}_{E_1} E_2=\nabla_{E_1}E_2+B\left(E_1,E_2\right)=0 \\
			&\hat{\nabla}_{E_2} E_2=\tilde{\nabla}_{E_2} E_2+\Phi=\nabla_{E_2}E_2+B\left(E_2,E_2\right)+\Phi=\frac{3}{4} \frac{E_1 f}{f} E_1+3 f E_3-c f^{3 / 2} E_4+\Phi \\
			&\hat{\nabla}_{E_1} E_3=\tilde{\nabla}_{E_1}E_3=-A_3\left(E_1\right)+\nabla^\perp_{E_1}E_3=f E_1 \\
			&\hat{\nabla}_{E_2} E_3=\tilde{\nabla}_{E_2} E_3=-A_3\left(E_2\right)+\nabla^\perp_{E_2}E_3=-3 f E_2 \\
			&\hat{\nabla}_{E_1} E_4=\tilde{\nabla}_{E_1} E_4=-A_4\left(E_1\right)+\nabla^\perp_{E_1}E_4=-c f^{3 / 2} E_1 \\
			&\hat{\nabla}_{E_2} E_4=\tilde{\nabla}_{E_2} E_4=-A_4\left(E_2\right)+\nabla^\perp_{E_2}E_4=c f^{3 / 2} E_2 \\
			&\hat{\nabla}_{E_1} \Phi=E_1 \\
			&\hat{\nabla}_{E_2} \Phi=E_2.
		\end{aligned}
		\right.
	\end{equation} 
	Our aim is to find extrinsic properties of $M$ and, finally, to infer a parametrization of $M$.
	
	Recall that, in \cite{SNMR2024}, where we studied the PNMC biconservative surfaces $M$ in $\mathbb{H}^4$, mainly from an intrinsic point of view, we preferred to choose a local chart built on the flow of $E_1$, since the integral curves of $E_1$ are geodesics of $M$. This choice proved to be useful. In the present paper, we look for a parametrization of the PNMC biconservative surfaces in $\mathbb{H}^4$, viewed as surfaces in $\mathbb{R}_1^5$. For this purpose we need to know the geometric properties of the integral curves of $E_1$ and $E_2$, thought of as curves in $\mathbb{R}_1^5$. As we will see, the integral curves of $E_1$ in $\mathbb{R}_1^5$ (and also in $\mathbb{H}^4$), are no longer geodesics and are much more complicated than the integral curves of $E_2$ in $\mathbb{R}_1^5$. The latter lie in two-dimensional affine subspaces while the integral curves of $E_1$ do not enjoy this simple property. More precisely, in Propositions \ref{integralcurveE2-case1}, \ref{integralcurveE2-case2}, and \ref{integralcurveE2-case3}, we will provide a complete classification of the integral curves of $E_2$, viewed as curves in $\mathbb{R}_1^5$, according to the type of the vector field $\hat{\nabla}_{E_2}E_2$. If $\hat{\nabla}_{E_2}E_2$ is a non-zero null, spacelike, or timelike vector field, then the corresponding integral curves of $E_2$ are, respectively, affine parabolas, circles, and hyperbolas. For this reason, we prefer to choose a local chart constructed along the flow of $E_2$ instead of $E_1$.
	
	Let $p_0 \in M$ be an arbitrarily fixed point of $M$ and let $\sigma=\sigma\left(u\right)$ be an integral curve of $E_1$ with $\sigma\left(0\right)=p_0$. Considering $\left\{\phi_v\right\}_{v \in \mathbb{R}}$ the flow of $E_2$ near the point $p_0$, we can define the following positively oriented local chart on $M$
	$$
	X^f\left(u,v\right)=\phi_v\left(\sigma\left(u\right)\right)=\phi_{\sigma\left(u\right)}\left(v\right).
	$$
	Then, we have
	\begin{equation*}
		\left\{
		\begin{array}{llll}
			X^f\left(u,0\right)=\sigma\left(u\right)\\
			X^f_u\left(u,0\right)=\sigma'\left(u\right)=E_1\left(u,0\right)\\
			X^f_v\left(u,v\right)=\phi'_{\sigma\left(u\right)}\left(v\right)=E_2\left(\phi_{\sigma\left(u\right)}\left(v\right)\right)=E_2\left(u,v\right),
		\end{array}
		\right.
	\end{equation*}
	for any $u$ and $v$. As the mean curvature function $f$ satisfies $E_2 f=0$, it follows that $f$ depends only on $u$ and $f'\left(u\right)>0$, for any $u$. Moreover, using \eqref{conditions-f}, it follows that 
	$$
	\left(E_1f\right)\left(u,v\right)=f'\left(u\right)\qquad \text{and} \qquad \left(E_1\left(E_1f\right)\right)\left(u,v\right)=f''\left(u\right), 
	$$
	for any $u$ and $v$. Then, since $\grad f=\left(E_1f\right)E_1$ and  
	\begin{equation*}
		\Delta f =-E_1\left(E_1f\right)+\frac{3}{4}\frac{\left(E_1f\right)^2}{f},
	\end{equation*}
	equation \eqref{second-order-invariant-f} becomes
	\begin{equation}\label{second-order-f-E_1}
		fE_1\left(E_1 f\right)-\frac{7}{4}\left(E_1f\right)^2+\frac{4}{3}f^2+4f^4+\frac{4}{3}c^2f^5=0,
	\end{equation}
	or, equivalently,
	\begin{equation}\label{second-order-chart-f}
		ff''-\frac{7}{4}\left(f'\right)^2+\frac{4}{3}f^2+4f^4+\frac{4}{3}c^2f^5=0.
	\end{equation}
	Using the same technique as in \cite{SNMR2024}, we achieve the first integral of the above second-order ordinary differential equation
	\begin{equation}\label{first-integral-f}
		f'=\frac{4}{3}f\sqrt{1+9Cf^{3/2}-9f^2-c^2f^3}>0,
	\end{equation}
	or, equivalently,
	\begin{equation}\label{first-integral-f-E1}
		E_1f=\frac{4}{3}f\sqrt{1+9Cf^{3/2}-9f^2-c^2f^3}>0,
	\end{equation}
	where $C$ is a constant of integration. 
	\begin{remark}
		We note that equation \eqref{second-order-chart-f} coincides with $\left(3.5\right)$ in \cite{SNMR2024} (we use different local charts but in both equations, $u$ denotes the arc-length parameter for the integral curves of $E_1$). The constant $C$ appearing in \eqref{first-integral-f} represents $C/8$ appearing in $\left(3.6\right)$ in \cite{SNMR2024}.
	\end{remark}
	We continue with the study of geometric properties of the integral curves of $E_1$, viewed as curves in $\mathbb{H}^4$. We denote such a curve by $\sigma=\sigma\left(u\right)$. We want to find the Frenet frame field for $\sigma$. For a homogeneous notation, we denote by $V_1$ the restriction of $E_1$ along $\sigma$, i.e., 
	\begin{equation}\label{V1}
		V_1=E_1,
	\end{equation}
	and then, the following real-valued functions and vector fields  
	\begin{equation*}\label{kappa1}
		\tilde{\kappa}_1=\left|\tilde{\nabla}_{V_1} V_1\right|=f \sqrt{1+c^2 f}>0,
	\end{equation*}
	\begin{equation}\label{V2}
		V_2=\frac{1}{\tilde{\kappa}_1} \tilde{\nabla}_{V_1} V_1=\frac{1}{\sqrt{1+c^2 f}}\left(-E_3+c \sqrt{f} E_4\right),
	\end{equation}
	$$
	\tilde{\kappa}_2=\left|\tilde{\nabla}_{V_1} V_2+\tilde{\kappa}_1 V_1\right|=\frac{|c| f^{\prime}}{2 \sqrt{f}\left(1+c^2 f\right)}>0
	$$
	and 
	\begin{equation}\label{V3}
		V_3=\frac{1}{\tilde{\kappa}_2}\left(\tilde{\nabla}_{V_1}V_2+\tilde{\kappa}_1 V_1\right)=\frac{\sgn\left(c\right)}{\sqrt{1+c^2 f}}\left(c\sqrt{f}E_3+E_4\right).
	\end{equation}
	With these notations, we get that the orthonormal frame fields $\left\{V_i\right\}_{i\in\overline{1,3}}$ satisfy
	$$
	\left\{\begin{array}{l}
		\tilde{\nabla}_{V_1} V_1=\tilde{\kappa}_1 V_2 \\
		\tilde{\nabla}_{V_1} V_2=-\tilde{\kappa}_1 V_1+\tilde{\kappa}_2 V_3 \\
		\tilde{\nabla}_{V_1} V_3=-\tilde{\kappa}_2 V_2
	\end{array}\right..
	$$
	Now, we define the vector field $V_4$ by the requiring that $\left\{V_i\right\}_{i\in\overline{1,4}}$ be a positively oriented orthonormal frame field in $\mathbb{H}^4$ along $\sigma$, and thus $\left\{V_i\right\}_{i\in\overline{1,4}}$ becomes the Frenet frame field for $\sigma$ with
	$$
	\tilde{\nabla}_{V_1} V_3=-\tilde{\kappa}_2 V_2+\tilde{\kappa}_3 V_4,\\
	$$
	where $\tilde{\kappa}_3=0$. As $\tilde{\kappa}_3=0$, it follows that $\sigma$ lies in a totally geodesic hypersurface  
	$\mathbb{H}^3$ of $\mathbb{H}^4$, $\mathbb{H}^3=\mathbb{H}^4 \cap \Pi$, where $\Pi$ is a Lorentzian hyperplane of $\mathbb{R}^5_1$ which contains the origin. We mention that the normal direction to $\Pi$ is given by a spacelike vector. More precisely, since
	$$
	\hat{\nabla}_{V_1} E_2=\hat{\nabla}_{E_1} E_2=0,
	$$
	i.e., $E_2$ is constant along $\sigma$, and since $E_2$ is also orthogonal to $\sigma$ along $\sigma$, it follows that $E_2$ is normal to the hyperplane $\Pi$.
	
	Moreover, using \eqref{first-integral-f}, we obtain
	\begin{equation*}\label{kappa2}
		\tilde{\kappa}_2=\frac{2|c|\sqrt{f}}{3 \left(1+c^2 f\right)} \sqrt{1+9 C f^{3 / 2}-9 f^2- c^2 f^3}.
	\end{equation*}
	Next, we show that the integral curves $\sigma$ of $E_1$ in $\mathbb{R}^5_1$ are far from being geodesics; more precisely, they cannot lie in any two-dimensional affine subspace of $\mathbb{R}^5_1$.
	
	First, let $F_1$ be the restriction of $E_1$ along $\sigma$. Then, we have
	$$
	\hat{\nabla}_{F_1}F_1=\tilde{\nabla}_{V_1}V_1+\Phi=\tilde{\kappa}_1 V_2+\Phi
	$$
	and
	$$
	\langle \hat{\nabla}_{F_1}F_1,\hat{\nabla}_{F_1}F_1 \rangle =\tilde{\kappa}_1^2-1\neq 0.%=f^2\left(1+c^2f\right)-1%
	$$
	We now introduce the real-valued function
	\begin{equation*}\label{kappa-hat}
		\hat{\kappa}_1=\sqrt{\left|\langle \hat{\nabla}_{F_1}F_1,\hat{\nabla}_{F_1}F_1 \rangle\right|}=\sqrt{\left|\tilde{\kappa}_1^2-1\right|}>0
	\end{equation*}
	and the vector field 
	$$
	F_2=\frac{1}{\hat{\kappa}_1}\hat{\nabla}_{F_1}F_1.
	$$
	A standard computation then yields
	$$
	\hat{\nabla}_{F_1}F_2=-\frac{\tilde{\kappa}_1^2-1}{\hat{\kappa}_1}V_1+\left(\frac{\tilde{\kappa}_1}{\hat{\kappa}_1}\right)'V_2+\frac{\tilde{\kappa}_1\tilde{\kappa}_2}{\hat{\kappa}_1}V_3+\left(\frac{1}{\hat{\kappa}_1}\right)'\Phi.
	$$
	Moreover, we can check that 
	$$
	\hat{\kappa}_2=\sqrt{\left|\langle \hat{\nabla}_{F_1}F_2+\hat{\kappa}_1\langle F_2,F_2\rangle F_1,\hat{\nabla}_{F_1}F_2 +\hat{\kappa}_1\langle F_2,F_2\rangle F_1 \rangle\right|}\neq 0.
	$$
	%It is different from zero since otherwise we would obtain that $f$ is constant.
	Therefore, $\sigma$ cannot lie in any two-dimensional affine subspace of $\mathbb{R}^5_1$.
	
	We finally conclude with the following proposition related to the integral curves of $E_1$.
	\begin{proposition} \label{Proposition 3.2}
		Let $\varphi:M^2 \rightarrow \mathbb{H}^4$ be a PNMC biconservative surface and consider $\sigma=\sigma\left(u\right)$ to be an integral curve of $E_1$, viewed as a curve in $\mathbb{H}^4$. Then, the following hold:
		\begin{itemize}
			\item[(i)] $E_2$ is constant along $\sigma$;
			\item[(ii)] $\sigma$ lies in a totally geodesic hypersurface $\mathbb{H}^3=\mathbb{H}^4 \cap \Pi$, where the hyperplane $\Pi$ contains the origin and is orthogonal to $E_2$;
			\item[(iii)] $\sigma$ does not lie in any two-dimensional affine subspace of $\mathbb{R}^5_1$;
			\item[(iv)] the curvature and the torsion of $\sigma$ are
			$$
			\kappa\left(u\right)=\tilde{\kappa}_1\left(u\right)=f\left(u\right) \sqrt{1+c^2 f\left(u\right)}
			$$
			and
			$$
			\tau\left(u\right)=\tilde{\kappa}_2\left(u\right)=\frac{2|c|\sqrt{f\left(u\right)}}{3 \left(1+c^2 f\left(u\right)\right)} \sqrt{1+9 C f^{3 / 2}\left(u\right)-9 f^2\left(u\right)- c^2 f^3\left(u\right)},
			$$
			where $f$ is the mean curvature function, $c\in\mathbb{R}^\ast$ and $C\in\mathbb{R}$.
		\end{itemize}
	\end{proposition}
	In the following, we study the geometric properties of the integral curves of $E_2$, viewed as curves in $\mathbb{R}^5_1$, by determining their Frenet frame field. 
	
	Using the expression of $ \hat{\nabla}_{E_2} E_2$ from \eqref{formulas Th.3.1}, and \eqref{first-integral-f-E1}, one obtains 
	$$
	\langle\hat{\nabla}_{E_2} E_2, \hat{\nabla}_{E_2} E_2\rangle=9C f^{3 / 2}.
	$$
	Denote by $\hat{\kappa}$ the curvature of the integral curves of $E_2$ in $\mathbb{R}_1^5$, i.e.,
	\begin{equation*}
		\hat{\kappa}=\sqrt{\left|\langle\hat{\nabla}_{E_2} E_2, \hat{\nabla}_{E_2} E_2\rangle\right|}=3\sqrt{|C|} f^{3/4}.
	\end{equation*}
	Further, since $C$ is an arbitrary real constant, we distinguish three cases: $C=0$, $C>0$, and $C<0$, which correspond, respectively, to $\hat{\nabla}_{E_2} E_2$ being a null, spacelike, and timelike vector field.
	
	\subsection{The case $C=0$}
	%%%%%%%%%%%%%%%%%%%%%%%%%%%%%%%%%%%%%%%%%%%%%%%%%%%%%%%%%%%%%%%%%%%%%%%%%%%%
	Clearly, the curvature $\hat{\kappa}$ of the integral curves of $E_2$ is zero and then, we define the vector field
	\begin{equation*}
		\xi= \hat{\nabla}_{E_2} E_2.
	\end{equation*}
	Using \eqref{formulas Th.3.1}, \eqref{second-order-f-E_1} and \eqref{first-integral-f-E1}, by some straightforward computations, we get $\langle \xi,\xi\rangle=0$,
	\begin{equation*}
		\hat{\nabla}_{E_1} \xi=\sqrt{1-9f^2-c^2f^3}\xi \qquad\text{and}\qquad
		\hat{\nabla}_{E_2} \xi=0.
	\end{equation*}
	Now, we are ready to find the local parametrization of $M^2$ in $\mathbb{R}_1^5$. This parametrization will rely on a solution $f$ of a second-order ODE and on a certain curve in $\mathbb{H}^3$, uniquely determined by $f$ and the condition that its position vector has to make a specific angle with a constant direction.
	
	First, let us consider an integral curve $\gamma=\gamma\left(v\right)$ of $E_2$, viewed in $\mathbb{R}^5_1$. Then, as $\hat{\nabla}_{E_2}\xi=0$, it follows that $\xi$ is a constant vector field in $\mathbb{R}^5_1$, along $\gamma$, i.e., $\xi\left(v\right)=\xi\left(0\right)$, for any $v$. Moreover, since
	$$
	\gamma''\left(v\right)=\hat{\nabla}_{\gamma'}\gamma'=\xi\left(v\right)=\xi\left(0\right),
	$$
	one obtains $\gamma'''\left(v\right)=0$, for any $v$ and, more precisely, 
	\begin{equation}\label{gamma polinomial in v}
		\gamma\left(v\right)=\overline{B}_0+v\overline{B}_1+\frac{1}{2}v^2 \overline{B}_2,
	\end{equation}
	where $\overline{B}_0$, $\overline{B}_1$ and $\overline{B}_2$ are constant vectors in $\mathbb{R}_1^5$ given by
	$$
	\overline{B}_0=\gamma\left(0\right),\qquad \overline{B}_1=\gamma'\left(0\right)=E_2\left(0\right), \qquad \overline{B}_2=\xi\left(0\right).
	$$
	Since $\langle\gamma',\gamma'\rangle=1$, it follows that
	$$
	\left\langle \overline{B}_1, \overline{B}_1\right\rangle=1 \qquad\text{and}\qquad \left\langle \overline{B}_2, \overline{B}_2\right\rangle=\left\langle \overline{B}_1, \overline{B}_2\right\rangle=0.
	$$
	Therefore, we can state the following properties of the integral curves of $E_2$.
	\begin{proposition} \label{integralcurveE2-case1}
		Let $\varphi:M^2 \to \mathbb{H}^4$ be a PNMC biconservative surface, and let $\gamma=\gamma\left(v\right)$ be an integral curve of $E_2$, viewed as a curve in $\mathbb{R}_1^5$. Assume that $\hat{\nabla}_{E_2}E_2$ is a non-zero null vector field. Then, the following hold:
		\begin{itemize}
			\item[(i)] the curvature $\hat{\kappa}$ of $\gamma$ is zero;
			\item[(ii)] the vector field $\xi=\hat{\nabla}_{E_2}E_2$ is constant along $\gamma$;
			\item[(iii)] the curve $\gamma$ is given by
			$$
			\gamma\left(v\right)=\overline{B}_0+v\overline{B}_1+\frac{1}{2}v^2\overline{B}_2,
			$$
			where 
			$$\overline{B}_0=\gamma\left(0\right),\qquad \overline{B}_1=E_2\left(0\right),\qquad  \overline{B}_2=\xi
			$$
			and
			$$
			\left\langle \overline{B}_1,\overline{B}_1\right\rangle=1, \qquad 
			\left\langle \overline{B}_2,\overline{B}_2\right\rangle=
			\left\langle \overline{B}_1,\overline{B}_2\right\rangle=0.
			$$
			In particular, $\gamma$ is a parabola contained in the lightlike affine plane
			$$
			\overline{B}_0+\operatorname{span}\{\overline{B}_1,\overline{B}_2\}\subset \mathbb{R}_1^5.
			$$
		\end{itemize}
	\end{proposition}
	Now, we are ready to obtain the parametrization of $M$ as a surface in $\mathbb{R}^5_1$. Let $p_0 \in M$ be an arbitrarily fixed point of $M$ and $\sigma=\sigma\left(u\right)$, $u\in I$, be an integral curve of $E_1$, viewed in $\mathbb{R}^5_1$, with $\sigma\left(0\right)=p_0$, where $I\subset\mathbb{R}$ is an open interval containing $0$. Consider $\left\{\phi_v\right\}_{v \in \mathbb{R}}$ the flow of $E_2$ near the point $p_0$. Then, for any $u \in I$ and for any $v \in \mathbb{R}$, the parametrization $\Phi=\Phi\left(u, v\right)$ of $M$ in $\mathbb{R}^5_1$ is given by
	\begin{equation*}
		\Phi\left(u, v\right)=\phi_{v}\left(\sigma\left(u\right)\right)=\overline{B}_0\left(u\right)+v\overline{B}_1\left(u\right)+\frac{1}{2}v^2 \overline{B}_2\left(u\right),
	\end{equation*}
	where the vector-valued functions $\overline{B}_0, \overline{B}_1, \overline{B}_2$, which are uniquely determined by the surface, are given by
	\begin{equation*}
		\overline{B}_0\left(u\right)=\sigma\left(u\right), \qquad \overline{B}_1\left(u\right)=E_2\left(u,0\right),\qquad \overline{B}_2\left(u\right)=\xi\left(u,0\right), 
	\end{equation*}
	and they satisfy
	\begin{equation*}
		\left\langle \overline{B}_1\left(u\right), \overline{B}_1\left(u\right)\right\rangle=1\qquad\text{and}\qquad\left\langle \overline{B}_2\left(u\right), \overline{B}_2\left(u\right)\right\rangle=\left\langle \overline{B}_1\left(u\right), \overline{B}_2\left(u\right)\right\rangle=0.
	\end{equation*}
	Clearly, $\Phi\left(u,0\right)=\sigma\left(u\right)$. Moreover, as $E_2$ is a constant vector in $\mathbb{R}^5_1$ along $\sigma$, one obtains that 
	\begin{equation}\label{expression-b1-our-case}
		\overline{B}_1\left(u\right)=\overline{b}_1, \qquad\forall u, 
	\end{equation} 
	where $\overline{b}_1$ is a constant vector of $\mathbb{R}^5_1$ and $\langle\overline{b}_1,\overline{b}_1\rangle=1$.
	Concerning the vector-valued function $\overline{B}_2=\overline{B}_2\left(u\right)$, we first note that since
	$$
	\hat{\nabla}_{E_1} \xi=\sqrt{1-9f^2-c^2f^3}\xi
	$$
	and 
	$$
	\overline{B}_2\left(u\right)=\xi\left(u,0\right),
	$$
	we have
	$$
	\overline{B}_2'\left(u\right)-\sqrt{1-9f^2\left(u\right)-c^2f^3\left(u\right)}\ \overline{B}_2\left(u\right)=0.
	$$
	Now, using \eqref{first-integral-f}, the above ODE becomes
	$$
	\overline{B}_2'\left(u\right)-\frac{3}{4}\frac{f'\left(u\right)}{f\left(u\right)}\overline{B}_2\left(u\right)=0,
	$$
	and thus, we may consider
	\begin{equation}\label{expression-b2-our-case}
		\overline{B}_2\left(u\right)=2f^{3/4}\left(u\right)\overline{b}_2,
	\end{equation}
	where $\overline{b}_2\in\mathbb{R}^5_1$ and $\langle \overline{b}_2,\overline{b}_2\rangle=0$.
	With these remarks, the parametrization can be expressed as
	$$
	\Phi\left(u, v\right)=\sigma\left(u\right)+v\overline{b}_1+v^2f^{3/4}\left(u\right) \overline{b}_2, 
	$$
	where
	$$
	\langle\overline{b}_1,\overline{b}_1\rangle=1\qquad\text{and}\qquad\langle\overline{b}_2,\overline{b}_2\rangle=\langle \overline{b}_1,\overline{b}_2\rangle=0.
	$$
	From Proposition \ref{Proposition 3.2} we know that $E_2$ is constant along $\sigma$, so 
	$$
	\langle\overline{b}_1,\sigma\left(u\right)\rangle=0.
	$$ 
	Then, since $\langle \Phi\left(u,v\right),\Phi\left(u,v\right)\rangle=-1$, one obtains
	\begin{equation}\label{additional-inner-product-condition}
	\langle \sigma\left(u\right),\overline{b}_2\rangle=-\frac{1}{2f^{3/4}\left(u\right)}.
	\end{equation}
	Before stating the conclusion of this part, we verify that the additional inner-product condition \eqref{additional-inner-product-condition} imposed on the directrix curve $\sigma$ is compatible with its prescribed curvature and torsion.
	\begin{proposition}
		Let $\mathbb{H}^3$ be the totally geodesic hypersurface of $\mathbb{H}^4$ determined by the constant vector $\overline{b}_1\in\mathbb{R}^5_1$ given by \eqref{expression-b1-our-case}. Then, the torsion of $\sigma$ satisfies \eqref{general-tau} in Lemma \ref{lemma-1}, and the constant vector $\overline{b}_2\in\mathbb{R}^5_1$ given by \eqref{expression-b2-our-case} coincides with the vector $\overline{b}_2$ given by \eqref{expression-b2} in Lemma \ref{lemma-2}, for
		$$
		\delta=0 \qquad\text{and}\qquad \Theta\left(u\right)=-\frac{1}{2f^{3/4}\left(u\right)}.
		$$
	\end{proposition}
	\begin{proof}
		By a straightforward computation we obtain that the curvature and torsion of $\sigma$, as a curve in $\mathbb{H}^3$, given in Proposition \ref{Proposition 3.2}, satisfies \eqref{general-tau} for $\delta$ and $\Theta$ given above.
		
		For the second part of the proof, from \eqref{expression-b2-our-case} and \eqref{formulas Th.3.1}, and using the expressions of the Frenet frame fields for the integral curves of $E_1$ given in \eqref{V1}, \eqref{V2} and \eqref{V3}, we obtain
		\begin{align}
			\overline{b}_2 & =\frac{1}{2f^{3/4}\left(u\right)}\left(\sigma\left(u\right)+\frac{3}{4}\frac{f'\left(u\right)}{f\left(u\right)}E_1\left(u,0\right)+3f\left(u\right)E_3\left(u,0\right)-cf^{3/2}\left(u\right)E_4\left(u,0\right)\right) \nonumber\\
			&= \frac{1}{2f^{3/4}}\sigma+\frac{3}{8}\frac{f'}{f^{7/4}}V_1-\frac{f^{1/4} \left(3+c^2f\right)}{2\sqrt{1+c^2f}}V_2+\frac{|c|f^{3/4}}{\sqrt{1+c^2f}}V_3 \nonumber\\
			&=-\Theta\sigma+\Theta'V_1
			+\frac{\Theta''-\Theta}{\kappa}V_2\
			+\frac{
				\Theta'''+\left(\kappa^2-1\right)\Theta'
				-\frac{\kappa'}{\kappa}\left(\Theta''-\Theta\right)
			}{\kappa\tau}V_3.
		\end{align}
	\end{proof}
	We are now ready to state the main result of this part.
	\begin{theorem}\label{theorem-C=0}
		Let $\varphi:M^2 \rightarrow \mathbb{H}^4$ be a PNMC biconservative immersion and denote $\Phi=\iota \circ \varphi: M \rightarrow \mathbb{R}_1^5$, where $\iota: \mathbb{H}^4 \rightarrow \mathbb{R}_1^5$ is the canonical inclusion. Identifying $M$ with its image, the surface $M$ can be locally parametrized as
		$$
		\Phi\left(u, v\right)=\sigma\left(u\right)+v\overline{b}_1+v^2f^{3/4}\left(u\right)\overline{b}_2, 
		$$
		where
		\begin{itemize}
			\item[(i)] $f=f\left(u\right)$ is a positive solution of the first-order ordinary differential equation 
			$$
			f'=\frac{4}{3}f\sqrt{1-9f^2-c^2f^3}>0,
			$$
			where $c$ is a real constant;
			\item[(ii)] $\overline{b}_1$ is a constant vector in $\mathbb{R}_1^5$ such that $\left\langle \overline{b}_1, \overline{b}_1\right\rangle=1$;
			\item[(iii)] $\sigma=\sigma\left(u\right)$ is a curve  in $\mathbb{R}_1^5$ parametrized by arc-length which lies in a totally geodesic hypersurface $\mathbb{H}^3=\mathbb{H}^4 \cap \Pi$; the hyperplane $\Pi$ contains the origin and is orthogonal to $\overline{b}_1$. Moreover, the curvature and torsion of $\sigma$, as a curve in $\mathbb{H}^3$, are given in Proposition \ref{Proposition 3.2} and $\overline{b}_2$ is a constant vector in $\mathbb{R}_1^5$  which satisfies
		\end{itemize}
		$$
		\langle \overline{b}_2, \overline{b}_2\rangle=\langle \overline{b}_1, \overline{b}_2\rangle=0 \qquad \text{and}\qquad \left\langle\sigma\left(u\right), \overline{b}_2\right\rangle=-\frac{1}{2 f^{3/4}\left(u\right)}.
		$$
	\end{theorem}
	\begin{remark}
		We can say that the constant vector $\overline{b}_1$ fixes $\mathbb{H}^3$, the curvature and torsion given in Proposition \ref{Proposition 3.2} determine $\sigma$, and $\sigma$ determines (uniquely) the constant vector $\overline{b}_2$.
	\end{remark}
	In the following we will give some geometric properties of the parametrization given in Theorem \ref{theorem-C=0}.
	
	First, we denote by
	\begin{align*}
		O\left(4,1\right)
		& =
		\left\{
		L\in GL\left(5,\mathbb{R}\right) \, |\,
		\langle L\overline{x},L\overline{y}\rangle
		=
		\langle \overline{x},\overline{y}\rangle
		\text{ for all } \overline{x},\overline{y}\in\mathbb{R}_1^5
		\right\}\\
		& = \left\{
		L\in GL\left(5,\mathbb{R}\right) \, |\,
		L^T JL=J
		\right\}
	\end{align*}
	the Lorentz group of $\mathbb{R}_1^5$, where $J=\operatorname{diag}\left(1,1,1,1,-1\right)$, and by
	$$
	O^+\left(4,1\right)
	=
	\left\{
	L\in O\left(4,1\right) \, |\,
	L\left(\mathbb{H}^4\right)=\mathbb{H}^4
	\right\}
	$$
	the isometry group of $\mathbb{H}^4$. Thus, the elements of $O^+\left(4,1\right)$ are precisely the isometries of $\mathbb{R}^5_1$ that preserve $\mathbb{H}^4$.
	\begin{proposition}\label{prop-C=0}
		Let $\left\{\overline{e}_1=\overline{b}_1,\overline{e}_2,\ldots,\overline{e}_5\right\}$ be a pseudo-orthonormal basis in $\mathbb{R}_1^5$. We may assume that $\overline{b}_2=\overline{e}_4+\overline{e}_5$.
		Then, the isometries of $\mathbb{H}^4$ preserving the surface $M$ are precisely the transformations $N_s$ and $\tilde{J}N_s$, with $s\in\mathbb{R}$, where
		$$
		N_s=
		\left(
		\begin{array}{ccccc}
			1&0&0&-s&s\\\\
			0&1&0&0&0\\\\
			0&0&1&0&0\\\\
			s&0&0&1-\dfrac{s^2}{2}&\dfrac{s^2}{2}\\\\
			s&0&0&-\dfrac{s^2}{2}&1+\dfrac{s^2}{2}
		\end{array}
		\right)
		$$
		and $\tilde{J}=\operatorname{diag}\left(-1,1,1,1,1\right)$. The transformations $N_s$ form a one-parameter group of null rotations with null axis $\overline{b}_2$. More precisely,
		$$
		N_s\Phi\left(u,v\right)
		=
		\Phi\left(
		u,
		v+\frac{s}{2f^{3/4}\left(u\right)}
		\right),
		$$
		and
		$$
		\left(\tilde{J} N_s\right)\Phi\left(u,v\right)
		=
		\Phi\left(
		u,
		-v-\frac{s}{2f^{3/4}\left(u\right)}
		\right).
		$$
		In particular, for any fixed $u$, the curve $v\mapsto\Phi\left(u,v\right)$ is an orbit of the one-parameter group $\left\{N_s\, |\,  s\in\mathbb{R}\right\}$.
	\end{proposition}
	\begin{proof}
		With respect to the above pseudo-orthonormal basis, the parametrization becomes	
		$$
		\Phi\left(u,v\right)
		=
		\sigma\left(u\right)
		+
		v\overline{e}_1
		+
		v^2f^{3/4}\left(u\right)
		\left(\overline{e}_4+\overline{e}_5\right).
		$$
		Furthermore,
		$$
		\langle\sigma,\overline{e}_1\rangle
		=
		0,
		\qquad
		\langle\sigma,\overline{e}_4+\overline{e}_5\rangle
		=
		-\frac{1}{2f^{3/4}}.
		$$
		Let $L\in O^+\left(4,1\right)$ be an isometry of $\mathbb{H}^4$ preserving $M$. Such an isometry descends from a linear isometry of $\mathbb{R}^5_1$ and it preserves the norm of the mean curvature vector field, thus we have
		$$
		f\left(L\left(p\right)\right)=f\left(p\right)
		$$
		for any $p\in M$. Let
		$$
		L\left(p\right)=L\Phi\left(u,v\right)=\Phi\left(\widetilde{u},\widetilde{v}\right).
		$$
		It follows that
		$$
		f\left(\widetilde{u}\right)=f\left(u\right).
		$$
		As $f'\left(u\right)>0$, the function $f$ is strictly increasing and hence injective on $I$. Therefore, $\widetilde{u}=u$.
		
		Thus, for any fixed $u$, the isometry $L$ preserves the parabola
		$v\mapsto\Phi\left(u,v\right)$, which lies in the affine plane $\sigma\left(u\right)+\Span\left\{\overline{e}_1,\overline{e}_4+\overline{e}_5\right\}$.
		
		We have	
		$$
		\Phi_v
		=
		\overline{e}_1
		+
		2vf^{3/4}
		\left(\overline{e}_4+\overline{e}_5\right).
		$$
		Since	
		$$
		\langle\overline{e}_1,\overline{e}_1\rangle=1,
		\qquad
		\langle\overline{e}_1,\overline{e}_4+\overline{e}_5\rangle=0,
		\qquad
		\langle
		\overline{e}_4+\overline{e}_5,
		\overline{e}_4+\overline{e}_5
		\rangle
		=
		0,
		$$
		we obtain
		$$
		\langle\Phi_v,\Phi_v\rangle=1.
		$$
		Thus, $v$ is an arc-length parameter for the curve
		$v\mapsto\Phi\left(u,v\right)$. Since $L$ preserves the curve corresponding to each fixed $u$, there exists a reparametrization $\psi_u$ of it such that
		$$
		L\Phi\left(u,v\right)
		=
		\Phi\left(u,\psi_u\left(v\right)\right).
		$$
		Differentiating with respect to $v$ and using that $L$ is an isometry of $\mathbb{H}^4$, we obtain
		$$
		1
		=
		\left\langle
		L\Phi_v\left(u,v\right),
		L\Phi_v\left(u,v\right)
		\right\rangle
		=
		\left(\psi_u'\left(v\right)\right)^2
		\left\langle
		\Phi_v\left(u,\psi_u\left(v\right)\right),
		\Phi_v\left(u,\psi_u\left(v\right)\right)
		\right\rangle
		=\left(\psi_u'\left(v\right)\right)^2.
		$$
		Therefore
		$$
		\psi_u\left(v\right)
		=
		\mu\left(u\right)v+a\left(u\right),
		$$
		where
		$$
		\mu\left(u\right)=\mu\in\{-1,1\},
		\qquad
		a\left(u\right)\in\mathbb{R}.
		$$
		Here, $a\left(u\right)$ represents the translation of the arc-length parameter $v$ or $-v$ along the parabola $v\mapsto\Phi\left(u,v\right)$ corresponding to the fixed value of $u$.
		
		Hence,
		\begin{equation}\label{L-phi(u,v)}
			L\Phi\left(u,v\right)
			=
			\Phi\left(u,\mu v+a\left(u\right)\right),
			\qquad
			\mu\in\{-1,1\},
		\end{equation}
		for any $u$ and $v$. 
		
		Next, we will find the expressions of $a\left(u\right)$ and of $L$. Fixing $u$ arbitrarily, we expand both sides of \eqref{L-phi(u,v)} and compare the coefficients of $v^2$, $v$, and the constant terms. It follows
		$$
		L\left(\overline{e}_4+\overline{e}_5\right)
		=
		\overline{e}_4+\overline{e}_5,
		$$
		$$
		L\overline{e}_1
		=
		\mu\overline{e}_1
		+
		2\mu f^{3/4}\left(u\right)a\left(u\right)
		\left(\overline{e}_4+\overline{e}_5\right),
		$$
		and
		$$
		L\sigma\left(u\right)
		=
		\sigma\left(u\right)
		+
		a\left(u\right)\overline{e}_1
		+
		f^{3/4}\left(u\right)a^2\left(u\right)
		\left(\overline{e}_4+\overline{e}_5\right).
		$$
		Since $\overline{e}_i$ and $L\overline{e}_1$ are independent of $u$, there exists a real constant $s$ such that
		$$
		2\mu f^{3/4}\left(u\right)a\left(u\right)=s,
		$$
		or, equivalently,
		$$
		a\left(u\right)
		=
		\frac{\mu s}{2f^{3/4}\left(u\right)}.
		$$
		Therefore,
		$$
		L\Phi\left(u,v\right)=\Phi\left(u, \mu\left(v+\frac{s}{2f^{3/4}\left(u\right)}\right)\right),
		$$
		\begin{equation}\label{Le1}
			L\overline{e}_1
			=
			\mu\overline{e}_1
			+
			s\left(\overline{e}_4+\overline{e}_5\right),
		\end{equation}
		and
		$$
		L\sigma\left(u\right)
		=
		\sigma\left(u\right)
		+
		\frac{\mu s}{2f\left(u\right)^{3/4}}\overline{e}_1
		+
		\frac{s^2}{4f\left(u\right)^{3/4}}
		\left(\overline{e}_4+\overline{e}_5\right).
		$$
		Thus, for any fixed $u$, the vector $\left(L-\Id\right)\sigma\left(u\right)\in\Span\left\{\overline{e}_1,\overline{e}_4+\overline{e}_5\right\}$. On the other hand, as $\sigma$ lies full in $\mathbb{H}^3$, which is determined by $\overline{e}_1^\perp$, we have
		$$
		\Span
		\left\{
		\sigma\left(u\right)\, |\, u\in I
		\right\}
		=
		\overline{e}_1^{\perp}.
		$$
		Consequently, we obtain
		$$
		\left(L-\Id\right)
		\left(\overline{e}_1^{\perp}\right)
		\subset
		\Span
		\left\{
		\overline{e}_1,
		\overline{e}_4+\overline{e}_5
		\right\}.
		$$
		In particular, we may write
		$$
		L\overline{e}_j
		=
		\overline{e}_j
		+
		a_j\overline{e}_1
		+
		c_j\left(\overline{e}_4+\overline{e}_5\right), \qquad j\in\{2,3,4,5\}.
		$$
		First, we determine $L\overline{e}_4$. Since $L$ preserves the Lorentzian scalar product, we have
		$$
		\langle
		L\overline{e}_1,
		L\overline{e}_4
		\rangle
		=
		\langle
		\overline{e}_1,
		\overline{e}_4
		\rangle
		=
		0.
		$$
		On the other hand, using the expression for $L\overline{e}_1$ given in \eqref{Le1}, we obtain $s+\mu a_4=0$, and therefore $a_4=-\mu s$. Moreover,
		$$
		\langle
		L\overline{e}_4,
		L\overline{e}_4
		\rangle
		=
		\langle
		\overline{e}_4,
		\overline{e}_4
		\rangle
		=
		1.
		$$
		This gives $1+s^2+2c_4=1$, and consequently $c_4=-s^2/2$. Therefore,
		$$
		L\overline{e}_4
		=
		-\mu s\,\overline{e}_1
		+
		\left(1-\frac{s^2}{2}\right)\overline{e}_4
		-
		\frac{s^2}{2}\overline{e}_5.
		$$
		Similarly, writing
		$$
		L\overline{e}_5
		=
		\overline{e}_5
		+
		a_5\overline{e}_1
		+
		c_5\left(\overline{e}_4+\overline{e}_5\right).
		$$
		and using $\langle L\overline{e}_1, L\overline{e}_5 \rangle=0$ and $\langle L\overline{e}_5,
		L\overline{e}_5\rangle=-1$, we obtain
		$$
		L\overline{e}_5
		=
		\mu s\,\overline{e}_1
		+
		\frac{s^2}{2}\overline{e}_4
		+
		\left(1+\frac{s^2}{2}\right)\overline{e}_5.
		$$
		Next, we consider
		$$
		L\overline{e}_2
		=
		\overline{e}_2
		+
		a_2\overline{e}_1
		+
		c_2\left(\overline{e}_4+\overline{e}_5\right).
		$$
		Since $\langle L\overline{e}_1, L\overline{e}_2 \rangle =0$ and $\langle L\overline{e}_4, L\overline{e}_2 \rangle = 0$, one gets
		$$
		L\overline{e}_2
		=
		\overline{e}_2.
		$$
		Finally, in the same way,
		$$
		L\overline{e}_3
		=
		\overline{e}_3.
		$$
		Therefore, with respect to the pseudo-orthonormal basis $
		\{
		\overline{e}_1,
		\overline{e}_2,
		\overline{e}_3,
		\overline{e}_4,
		\overline{e}_5
		\}$,
		the matrix of $L$ is
		$$
		L_{\mu,s}
		=
		\begin{pmatrix}
			\mu&0&0&-\mu s&\mu s\\\\
			0&1&0&0&0\\\\
			0&0&1&0&0\\\\
			s&0&0&1-\frac{s^2}{2}&\frac{s^2}{2}\\\\
			s&0&0&-\frac{s^2}{2}&1+\frac{s^2}{2}
		\end{pmatrix}.
		$$
		If $\mu=1$, then $L_{1,s}=N_s$. The transformation $N_s$ fixes the null vector
		$\overline{b}_2=\overline{e}_4+\overline{e}_5$ and for any $s\neq 0$, we have 
		$$
		\left(N_s-\Id\right)\left(\overline{b}_2^\perp\right)= \mathbb{R}\overline{b}_2.
		$$
		Indeed,
		$$
		N_s\overline{e}_1
		=
		\overline{e}_1+s\overline{b}_2,
		\qquad
		N_s\overline{e}_2=\overline{e}_2,
		\qquad
		N_s\overline{e}_3=\overline{e}_3,
		\qquad
		N_s\overline{b}_2=\overline{b}_2.
		$$
		Hence, $N_s$ is a null rotation with null axis $\overline{b}_2$.
		
		If $\mu=-1$, let $\tilde{J}=\operatorname{diag}\left(-1,1,1,1,1\right)$. Then, 
		$$
		L_{-1,s}=\tilde{J} N_s=N_{-s} \tilde{J}.
		$$
		We mention that the transformation $\tilde{J}$ is the Lorentzian reflection with respect to $\overline{e}_1^{\perp}$.
		
		Conversely, the matrices $N_s$ and $\tilde{J}N_s$ belong to $O^+\left(4,1\right)$. Moreover, we can check directly that
		$$
		N_s\Phi\left(u,v\right)
		=
		\Phi\left(
		u,
		v+\frac{s}{2f\left(u\right)^{3/4}}
		\right)
		$$	
		and
		$$
		\left(\tilde{J} N_s\right)\Phi\left(u,v\right)
		=
		\Phi\left(
		u,
		-v-\frac{s}{2f\left(u\right)^{3/4}}
		\right).
		$$
		Hence, both $N_s$ and $\tilde{J}N_s$ preserve $M$. Therefore, these are all the isometries of $\mathbb{H}^4$ preserving the surface.
		
		Finally, we note that 
		$$
		N_sN_t=N_{s+t},
		\qquad
		N_0=I_5,
		$$
		so the transformations $N_s$ form a one-parameter group. For $s\neq0$, $N_s$ is a nontrivial null rotation fixing the null direction $\mathbb{R}\overline{b}_2$. Moreover, we can verify that $N_s$ has no fixed point in $\mathbb{H}^4$. Hence, its restriction to $\mathbb{H}^4$ is a parabolic isometry with a unique fixed point at infinity, corresponding to the null direction $\mathbb{R}\overline{b}_2$.
		
		Furthermore, for any fixed $u$ and any $v_0\in\mathbb{R}$,
		$$
		\left\{
		N_s\Phi\left(u,v_0\right)
		\,\middle|\,
		s\in\mathbb{R}
		\right\}
		=
		\left\{
		\Phi\left(u,v\right)
		\,\middle|\,
		v\in\mathbb{R}
		\right\}.
		$$
		Thus, the curve $v\mapsto\Phi\left(u,v\right)$ is an orbit of the one-parameter group $\{N_s\, |\, s\in\mathbb{R}\}$.
	\end{proof}
	
	\subsection{The case $C\neq 0$}
	%%%%%%%%%%%%%%%%%%%%%%%%%%%%%%%%%%%%%%%%%%%%%%%%%%%%%%%%%%%%%%%%%%%%%%%
	
	Clearly, the curvature $\hat{\kappa}$ of the integral curves of $E_2$ is positive and we can define the vector field
	\begin{equation*}\label{definition of xi}
		\xi=\frac{1}{\hat{\kappa}} \hat{\nabla}_{E_2} E_2=\frac{1}{3\sqrt{|C|} f^{3 / 4}}\hat{\nabla}_{E_2} E_2.
	\end{equation*}
	As $E_2f=0$ and $E_2\left(E_1 f\right)=0$, by standard computations, using also \eqref{formulas Th.3.1} and \eqref{first-integral-f-E1}, we obtain
	$$
	\hat{\nabla}_{E_2} \hat{\nabla}_{E_2} E_2=-9Cf^{3/2}E_2=-\sgn\left(C\right)\hat{\kappa}^2E_2.
	$$
	These relations lead to 
	$$
	\frac{E_1\hat{\kappa}}{\hat{\kappa}}=\frac{3}{4}\frac{E_1 f}{f}, \qquad E_2\hat{\kappa}=0,
	$$ 
	and 
	\begin{equation}\label{hat-nabla-xi}
		\hat{\nabla}_{E_1}\xi=0, \qquad \hat{\nabla}_{E_2}\xi=\frac{1}{\hat{\kappa}}\hat{\nabla}_{E_2} \hat{\nabla}_{E_2} E_2=-\sgn\left(C\right)\hat{\kappa}E_2.
	\end{equation}
	\begin{remark}
		Since $E_2\hat{\kappa}=0$,
		$$
		\hat{\nabla}_{E_2} E_2=\hat{\kappa}\xi \qquad \text{and} \qquad \hat{\nabla}_{E_2}\xi=-\sgn\left(C\right)\hat{\kappa}E_2,
		$$
		it follows that the integral curves of $E_2$ have constant curvature and they are contained in the two-dimensional affine subspace of $\mathbb{R}_1^5$ generated by $\xi\left(0\right)$ and $E_2\left(0\right)$.
	\end{remark}

	\subsubsection{The subcase $C>0$}
		%%%%%%%%%%%%%%%%%%%%%%%%%%%%%%%%%%%%%%%%%%%%%%%%%%%%%%%%%%%%%
		
	First, let us consider an integral curve $\gamma=\gamma\left(v\right)$ of $E_2$, viewed in $\mathbb{R}^5_1$. Then, as $\hat{\nabla}_{E_2}E_2=\hat{\kappa}\xi$ and $\hat{\nabla}_{E_2}\xi=-\hat{\kappa}E_2$, it follows that, along $\gamma$, we have
	\begin{equation*}
		\left\{
		\begin{array}{l}
			\gamma''\left(v\right)=\hat{\kappa}\xi\left(v\right)\\
			\xi'\left(v\right)=-\hat{\kappa}\gamma'\left(v\right)
		\end{array}
		\right..
	\end{equation*}
	Differentiating the first equation of the above system, one obtains the ODE
	$$
	\gamma'''\left(v\right)+\hat{\kappa}^2\gamma'\left(v\right)=0.
	$$
	A standard computation leads to
	$$
	\gamma\left(v\right)=\overline{A}_0+\frac{\sin\left(\hat{\kappa}v\right)}{\hat{\kappa}}\overline{A}_1-\frac{\cos\left(\hat{\kappa}v\right)}{\hat{\kappa}}\overline{A}_2,
	$$
	where $\overline{A}_0$, $\overline{A}_1$ and $\overline{A}_2$ are constant vectors in $\mathbb{R}_1^5$ 
	given by
	$$
	\overline{A}_0=\gamma\left(0\right)+\frac{1}{\hat{\kappa}}\xi\left(0\right),\qquad \overline{A}_1=\gamma'\left(0\right)=E_2\left(0\right), \qquad \overline{A}_2=\xi\left(0\right).
	$$
	Since $\langle\gamma',\gamma'\rangle=1$, it follows that
	$$
	\left\langle \overline{A}_1, \overline{A}_2\right\rangle=0 \qquad\text{and}\qquad \left\langle \overline{A}_1, \overline{A}_1\right\rangle=\left\langle \overline{A}_2, \overline{A}_2\right\rangle=1.
	$$
	If we denote by 
	$$
	\overline{B}_0=\overline{A}_0, \qquad \overline{B}_1=\frac{1}{\hat{\kappa}}\overline{A}_1 \qquad \text{and}\qquad \overline{B}_2=-\frac{1}{\hat{\kappa}} \overline{A}_2,  
	$$
	we rewrite
	\begin{equation*}
		\gamma\left(v\right)=\overline{B}_0+\sin\left(\hat{\kappa} v\right)\overline{B}_1+\cos\left(\hat{\kappa} v\right)\overline{B}_2,
	\end{equation*}
	where $\overline{B}_i \in \mathbb{R}_1^5$, $i=\overline{0,2}$ satisfy
	$$
	\overline{B}_0=\gamma\left(0\right)+\frac{1}{\hat{\kappa}}\xi\left(0\right),\qquad \overline{B}_1=\frac{1}{\hat{\kappa}}\gamma'\left(0\right)=\frac{1}{\hat{\kappa}}E_2\left(0\right), \qquad \overline{B}_2=-\frac{1}{\hat{\kappa}}\xi\left(0\right).
	$$
	Moreover, we have
	$$
	\left\langle \overline{B}_1,\overline{B}_2\right\rangle=0 \qquad\text{and}\qquad \left\langle \overline{B}_1,\overline{B}_1\right\rangle=\left\langle \overline{B}_2,\overline{B}_2\right\rangle=
	\frac{1}{\hat{\kappa}^2}.
	$$
	Thus, we can state the following properties of the integral curves of $E_2$.
	\begin{proposition}\label{integralcurveE2-case2}
		Let $\varphi:M^2 \to \mathbb{H}^4$ be a PNMC biconservative surface, and let $\gamma=\gamma\left(v\right)$ be an integral curve of $E_2$, viewed as a curve in $\mathbb{R}_1^5$. Assume that $\hat{\nabla}_{E_2}E_2$ is a non-zero spacelike vector field. Then, the following hold:
		\begin{itemize}
			\item[(i)] the curvature $\hat{\kappa}$ of $\gamma$ is a positive constant along $\gamma$;
			\item[(ii)] the vector field $\xi$ defined by $\xi=\left(1/\hat{\kappa}\right) \hat{\nabla}_{E_2} E_2$ satisfies $\hat{\nabla}_{E_2}\xi=-\hat{\kappa}E_2$;
			\item[(iii)] the curve $\gamma$ is given by
			$$
			\gamma\left(v\right)=\overline{B}_0+\sin\left(\hat{\kappa}v\right)\overline{B}_1+\cos\left(\hat{\kappa}v\right)\overline{B}_2,
			$$
			where
			$$
			\overline{B}_0=\gamma\left(0\right)+\frac{1}{\hat{\kappa}}\xi\left(0\right),\qquad
			\overline{B}_1=\frac{1}{\hat{\kappa}}E_2\left(0\right),\qquad
			\overline{B}_2=-\frac{1}{\hat{\kappa}}\xi\left(0\right)
			$$
			and
			$$
			\left\langle \overline{B}_1,\overline{B}_2\right\rangle=0,\qquad
			\left\langle \overline{B}_1,\overline{B}_1\right\rangle=
			\left\langle \overline{B}_2,\overline{B}_2\right\rangle=\frac{1}{\hat{\kappa}^2}.
			$$
			In particular, $\gamma$ is a circle of radius $1/\hat{\kappa}$ centered at $\overline{B}_0$ in the Euclidean affine plane
			$$
			\overline{B}_0+\operatorname{span}\{\overline{B}_1,\overline{B}_2\}\subset\mathbb{R}^5_1.
			$$
		\end{itemize}
	\end{proposition}
	To obtain the parametrization of $M$ as a surface in $\mathbb{R}^5_1$, let $p_0\in M$ be an arbitrarily fixed point of $M$ and $\sigma=\sigma\left(u\right)$ be an integral curve of $E_1$ with $\sigma\left(0\right)=p_0$. Consider $\left\{\phi_v\right\}_{v \in \mathbb{R}}$ the flow of $E_2$ near the point $p_0$. Then, for any $u \in\left(-\varepsilon, \varepsilon\right)$ and for any $v \in \mathbb{R}$, the parametrization $\Phi=\Phi\left(u, v\right)$ of $M$ in $\mathbb{R}^5_1$ is given by 
	$$
	\Phi\left(u, v\right)=\phi_{v}\left(\sigma\left(u\right)\right)=\overline{B}_0\left(u\right)+\sin\left(\hat{\kappa}\left(u\right)v\right)\overline{B}_1\left(u\right)+\cos \left(\hat{\kappa}\left(u\right)v\right)\overline{B}_2\left(u\right),
	$$
	where the vector-valued functions $\overline{B}_0, \overline{B}_1, \overline{B}_2$, which are uniquely determined by the surface, are given by
	\begin{equation*}
		\overline{B}_0\left(u\right)=\sigma\left(u\right)+\frac{1}{\hat{\kappa}\left(u\right)}\xi\left(u,0\right), \qquad \overline{B}_1\left(u\right)=\frac{1}{\hat{\kappa}\left(u\right)}E_2\left(u,0\right),\qquad \overline{B}_2\left(u\right)=-\frac{1}{\hat{\kappa}\left(u\right)}\xi\left(u,0\right),
	\end{equation*}
	and they satisfy
	\begin{equation*}
		\left\langle \overline{B}_1\left(u\right), \overline{B}_2\left(u\right)\right\rangle=0\qquad\text{and}\qquad\left\langle \overline{B}_1\left(u\right), \overline{B}_1\left(u\right)\right\rangle=\left\langle \overline{B}_2\left(u\right), \overline{B}_2\left(u\right)\right\rangle=\frac{1}{\hat{\kappa}^2\left(u\right)}.
	\end{equation*}
	Let us consider the vector-valued functions in $\mathbb{R}^5_1$ given by
	\begin{equation}\label{expressions-bi-second-case}
		\overline{b}_i\left(u\right)=\hat{\kappa}\left(u\right)\overline{B}_i\left(u\right),\qquad i=1,2.
	\end{equation}
	Then,
	$$
	\overline{b}_1\left(u\right)=E_2\left(u,0\right), \qquad \overline{b}_2\left(u\right)=-\xi\left(u,0\right),
	$$
	and the parametrization can be expressed as
	$$
	\Phi\left(u, v\right)=\sigma\left(u\right)+\frac{1}{\hat{\kappa}\left(u\right)}\left(\sin\left(\hat{\kappa}\left(u\right)v\right)\overline{b}_1\left(u\right)+\left(\cos \left(\hat{\kappa}\left(u\right)v\right)-1\right)\overline{b}_2\left(u\right)\right).
	$$
	Clearly, $\Phi\left(u,0\right)=\sigma\left(u\right)$. Moreover, by Proposition \ref{Proposition 3.2} and the first relation in \eqref{hat-nabla-xi}, it follows that $E_2$ and $\xi$  are constant vector fields in $\mathbb{R}^5_1$ along $\sigma$. Therefore, the vector-valued functions $\overline{b}_i=\overline{b}_i\left(u\right)$ are constant in $\mathbb{R}^5_1$; hence, they can be identified with the constant vectors $\overline{b}_i\in\mathbb{R}^5_1$ satisfying
	$$
	\langle\overline{b}_1,\overline{b}_2\rangle=0\qquad \text{and}\qquad \langle\overline{b}_1,\overline{b}_1\rangle=\langle\overline{b}_2,\overline{b}_2\rangle=1.
	$$
	In order to get a simpler expression of $\Phi$ we can consider the following change of coordinates
	$\left(u, v\right) \rightarrow\left(u, t=\hat{\kappa}\left(u\right) v\right)$. With respect to these new local coordinates, the parametrization $\Phi$ can be expressed as
	$$
	\Phi\left(u, t\right)=\sigma\left(u\right)+\frac{1}{\hat{\kappa}\left(u\right)}\left(\sin t\ \overline{b}_1+\left(\cos t-1\right)\overline{b}_2\right).
	$$
	From Proposition \ref{Proposition 3.2}, we know that $E_2$ is constant along $\sigma$, so 
	$$
	\langle\overline{b}_1,\sigma\left(u\right)\rangle=0.
	$$ 
	Then, since $\langle \Phi\left(u,v\right),\Phi\left(u,v\right)\rangle=-1$, one obtains
	$$
	\langle \sigma\left(u\right),\overline{b}_2\rangle=\frac{1}{\hat{\kappa}\left(u\right)}.
	$$
	Before stating the conclusion of this second part, we state the following result that can be obtained by straightforward computations.
	\begin{proposition}
		Let $\mathbb{H}^3$ be the totally geodesic hypersurface of $\mathbb{H}^4$ determined by the constant vector $\overline{b}_1\in\mathbb{R}^5_1$ given by \eqref{expressions-bi-second-case}. Then, the torsion of $\sigma$ satisfies \eqref{general-tau} in Lemma \ref{lemma-1}, and the constant vector $\overline{b}_2\in\mathbb{R}^5_1$ given by \eqref{expressions-bi-second-case} coincides with the vector $\overline{b}_2$ given by \eqref{expression-b2} in Lemma \ref{lemma-2}, for
		$$
		\delta=1 \qquad\text{and}\qquad \Theta\left(u\right)=\frac{1}{3\sqrt{C}f^{3/4}\left(u\right)}.
		$$
	\end{proposition}
	We are now ready to state the main result of this second case.
	\begin{theorem}\label{theorem-C>0}
		Let $\varphi:M^2\rightarrow \mathbb{H}^4$ be a PNMC biconservative immersion and denote $\Phi=\iota \circ \varphi: M \rightarrow \mathbb{R}_1^5$, where $\iota: \mathbb{H}^4 \rightarrow \mathbb{R}_1^5$ is the canonical inclusion. Identifying $M$ with its image, the surface $M$ can be locally parametrized as
		$$
		\Phi\left(u, t\right)=\sigma\left(u\right)+\frac{1}{\hat{\kappa}\left(u\right)}\left(\sin t\ \overline{b}_1+\left(\cos t-1\right)\overline{b}_2\right),
		$$
		where
		\begin{itemize}
			\item[(i)] $f=f\left(u\right)$ is a positive solution of the first-order ordinary differential equation 
			$$
			f'=\frac{4}{3}f\sqrt{1+9Cf^{3/2}-9f^2-c^2f^3}>0,
			$$
			where $c$ is a real constant and $C$ is a positive real constant;
			\item[(ii)]$\hat{\kappa}=\hat{\kappa}\left(u\right)$ is the positive function given by
			$$
			\hat{\kappa}\left(u\right)=3\sqrt{C}f^{3/4}\left(u\right);
			$$
			\item[(iii)] $\overline{b}_1$ is a constant vector in $\mathbb{R}_1^5$ such that $\langle\overline{b}_1,\overline{b}_1\rangle=1$;
			\item[(iv)] $\sigma=\sigma\left(u\right)$ is a curve in $\mathbb{R}_1^5$ parametrized by arc-length which lies in a totally geodesic hypersurface $\mathbb{H}^3=\mathbb{H}^4 \cap \Pi$; the hyperplane $\Pi$ contains the origin and is orthogonal to $\overline{b}_1$. Moreover, the curvature and torsion of $\sigma$, as a curve in $\mathbb{H}^3$, are given in Proposition \ref{Proposition 3.2} and $\overline{b}_2$ is a constant vector in $\mathbb{R}_1^5$ which satisfies
		\end{itemize}
		$$
		\langle\overline{b}_2,\overline{b}_1\rangle=0, \qquad \langle\overline{b}_2,\overline{b}_2\rangle=1 \qquad \text{and}\qquad  \langle\sigma\left(u\right), \overline{b}_2\rangle=\frac{1}{\hat{\kappa}\left(u\right)}.
		$$
	\end{theorem}
	Further, we give some geometric properties of the parametrization given in Theorem \ref{theorem-C>0}.
	\begin{proposition}
		Let $\{\overline{e}_1=\overline{b}_1,\overline{e}_2=\overline{b}_2,\overline{e}_3,\ldots,\overline{e}_5\}$ be a pseudo-orthonormal basis in $\mathbb{R}_1^5$. Then, the isometries of $\mathbb{H}^4$ preserving the surface $M$ are precisely the transformations $R_s$ and $\tilde{J}R_s$, with $s\in\mathbb{R}$, where
		$$
		R_s=
		\left(
		\begin{array}{ccccc}
			\cos s&-\sin s&0&0&0\\
			\sin s&\cos s&0&0&0\\
			0&0&1&0&0\\
			0&0&0&1&0\\
			0&0&0&0&1
		\end{array}
		\right)
		$$
		and $\tilde{J}=\operatorname{diag}\left(-1,1,1,1,1\right)$. The transformations $R_s$ form a one-parameter group of rotations in the spacelike plane $\Span\{\overline{e}_1,\overline{e}_2\}$. More precisely,
		$$
		R_s\Phi\left(u,t\right)
		=
		\Phi\left(u,t-s\right),
		$$
		and
		$$
		\left(\tilde{J} R_s\right)\Phi\left(u,t\right)
		=
		\Phi\left(u,-t+s\right).
		$$
		In particular, for any fixed $u$, the curve $t\mapsto\Phi\left(u,t\right)$ is an orbit of the one-parameter group $\{R_s\, | \, s\in\mathbb{R}\}$.
	\end{proposition}
	\begin{proof}
		With respect to the above pseudo-orthonormal basis, the parametrization becomes
		$$
		\Phi\left(u,t\right)
		=
		\sigma\left(u\right)
		+
		\frac{1}{\hat{\kappa}\left(u\right)}
		\left(
		\sin t\,\overline{e}_1
		+
		\left(\cos t-1\right)\overline{e}_2
		\right).
		$$
		Furthermore, we have	
		$$
		\langle\sigma,\overline{e}_1\rangle=0,
		\qquad
		\langle\sigma,\overline{e}_2\rangle
		=
		\frac{1}{\hat{\kappa}\left(u\right)}.
		$$
		If we denote
		\begin{equation}\label{alpha-sigma-C>0}
			\alpha\left(u\right)=\sigma\left(u\right)
			-
			\frac{1}{\hat{\kappa}\left(u\right)}\overline{e}_2,
		\end{equation}
		then 
		$$
		\alpha\left(u\right)
		\in
		\Span
		\{\overline{e}_3,\overline{e}_4,\overline{e}_5\},
		$$
		and the parametrization can be written as
		$$
		\Phi\left(u,t\right)
		=
		\alpha\left(u\right)
		+
		\frac{1}{\hat{\kappa}\left(u\right)}
		\left(
		\sin t\,\overline{e}_1
		+
		\cos t\,\overline{e}_2
		\right).
		$$
		Thus, for any fixed $u$, the curve $t\mapsto\Phi\left(u,t\right)$ is a circle with center $\alpha\left(u\right)$ and radius $1/\hat{\kappa}\left(u\right)$, lying in the affine plane $\alpha\left(u\right)+\Span\{\overline{e}_1,\overline{e}_2\}$.
		
		Let $L\in O^+\left(4,1\right)$ be an isometry of $\mathbb{H}^4$ preserving $M$. First, we note that, using the same argument as in the proof of Proposition \ref{prop-C=0}, $L$ preserves each circle $t\mapsto\Phi\left(u,t\right)$ and, in particular, its center. Hence,
		\begin{equation}\label{eq-L-dif-C=0}
			L\Phi\left(u,t\right)
			=
			\Phi\left(u,\psi_u\left(t\right)\right),
		\end{equation}
		where $\psi_u$ represents a reparametrization of the circle $t\mapsto\Phi\left(u,t\right)$, for an arbitrary fixed $u$.  
		Moreover, we can see that 
		$$
		\Phi_t=\frac{1}{\hat{\kappa}\left(u\right)}\left(\cos t\,\overline{e}_1-\sin t\, \overline{e}_2\right)
		$$	
		and 
		$$
		\langle\Phi_t,\Phi_t\rangle=\frac{1}{\hat{\kappa}^2\left(u\right)}.
		$$
		Differentiating \eqref{eq-L-dif-C=0} with respect to $t$ and using that $L$ is an isometry of $\mathbb{H}^4$, we obtain
		$$
		\frac{1}{\hat{\kappa}^2\left(u\right)}
		=
		\left\langle
		L\Phi_t\left(u,t\right),
		L\Phi_t\left(u,t\right)
		\right\rangle
		=
		\left(\psi_u'\left(t\right)\right)^2
		\left\langle
		\Phi_t\left(u,\psi_u\left(t\right)\right),
		\Phi_t\left(u,\psi_u\left(t\right)\right)
		\right\rangle
		=\left(\psi_u'\left(t\right)\right)^2\frac{1}{\hat{\kappa}^2\left(u\right)}.
		$$
		Therefore,
		$$
		\psi_u\left(t\right)
		=
		\mu\left(u\right)t+a\left(u\right),
		$$
		where
		$$
		\mu\left(u\right)=\mu\in\{-1,1\},
		\qquad
		a\left(u\right)\in\mathbb{R}.
		$$
		Hence,
		\begin{equation}\label{L-phi(u,t)-C>0}
			L\Phi\left(u,t\right)
			=
			\Phi\left(u,\mu t+a\left(u\right)\right),
			\qquad
			\mu\in\{-1,1\}.
		\end{equation}
		Next, we will determine $a\left(u\right)$ and $L$. Expanding both sides of \eqref{L-phi(u,t)-C>0}, we obtain
		\begin{equation}\label{extended-L-phi-C>0}
			L\alpha\left(u\right)
			+
			\frac{1}{\hat{\kappa}\left(u\right)}
			\left(
			\sin t\,L\overline{e}_1
			+
			\cos t\,L\overline{e}_2
			\right)
			=
			\alpha\left(u\right)
			+
			\frac{1}{\hat{\kappa}\left(u\right)}
			\left(
			\sin\left(\mu t+a\left(u\right)\right)\,\overline{e}_1
			+
			\cos\left(\mu t+a\left(u\right)\right)\,\overline{e}_2
			\right).
		\end{equation}
		Since $\mu\in\{-1,1\}$, we have
		$$
		\sin\left(\mu t+a\left(u\right)\right)
		=
		\mu\cos a\left(u\right)\sin t
		+
		\sin a\left(u\right)\cos t
		$$
		and
		$$
		\cos\left(\mu t+a\left(u\right)\right)
		=
		-\mu\sin a\left(u\right)\sin t
		+
		\cos a\left(u\right)\cos t.
		$$
		Therefore, comparing the coefficients of $\sin t$, $\cos t$, and the constant terms in \eqref{extended-L-phi-C>0}, we obtain
		$$
		L\overline{e}_1
		=\mu\cos a\left(u\right)\,\overline{e}_1-\mu\sin a\left(u\right)\,\overline{e}_2,
		$$
		$$
		L\overline{e}_2
		=\sin a\left(u\right)\,\overline{e}_1+\cos a\left(u\right)\,\overline{e}_2,
		$$
		and
		$$
		L\alpha\left(u\right)=\alpha\left(u\right).
		$$
		Since $\overline{e}_i$ and $L\overline{e}_i$ are independent of $u$, both $\sin a\left(u\right)$ and $\cos a\left(u\right)$ are independent of $u$. Hence, there exists a real constant $s$ such that
		$$
		a\left(u\right)=-\mu s.
		$$
		Consequently,
		$$
		L\Phi\left(u,t\right)
		=
		\Phi\left(u,\mu\left(t-s\right)\right),
		$$
		\begin{equation}\label{Le1-C>0}
			L\overline{e}_1
			=	
			\mu\cos s\,\overline{e}_1
			+
			\sin s\,\overline{e}_2,
		\end{equation}
		and
		\begin{equation}\label{Le2-C>0}
			L\overline{e}_2
			=
			-\mu\sin s\,\overline{e}_1
			+
			\cos s\,\overline{e}_2.
		\end{equation}
		Hence,
		$L\left(
		\Span
		\{\overline{e}_1,\overline{e}_2\}
		\right)
		=
		\Span
		\{\overline{e}_1,\overline{e}_2\}$.
		
		Since $L$ preserves the Lorentzian scalar product, it also preserves the Lorentz-orthogonal complement of the plane $\Span
		\{\overline{e}_1,\overline{e}_2\}$. Therefore,
		$$
		L\left(
		\Span
		\{\overline{e}_3,\overline{e}_4,\overline{e}_5\}
		\right)
		=
		\Span
		\{\overline{e}_3,\overline{e}_4,\overline{e}_5\}.
		$$
		We now prove that
		$$
		\Span 
		\left\{\alpha\left(u\right)\, |\,  u\in I \right\}
		=
		\Span
		\{\overline{e}_3,\overline{e}_4,\overline{e}_5\}.
		$$
		First, we recall that as $\sigma$ lies full in $\mathbb{H}^3$, which is determined by $\overline{e}_1^\perp$, we have 
		$$
		\Span \{\sigma\left(u\right)\, |\,  u\in I\}= \overline{e}_1^{\perp}=\Span\left\{\overline{e}_2,\overline{e}_3,\overline{e}_4,\overline{e}_5\right\}.
		$$ 
		As $\alpha\left(u\right)$ is given in \eqref{alpha-sigma-C>0} and the $\overline{e}_2$-component of $\alpha\left(u\right)$ vanishes we obtain that 
		$$
		\Span 
		\left\{\alpha\left(u\right)\, |\,  u\in I \right\}
		=
		\Span
		\{\overline{e}_3,\overline{e}_4,\overline{e}_5\}.
		$$
		Now, since $L$ fixes every vector $\alpha\left(u\right)$ we obtain
		$$
		L\overline{e}_3=\overline{e}_3,
		\qquad
		L\overline{e}_4=\overline{e}_4,
		\qquad
		L\overline{e}_5=\overline{e}_5.
		$$
		Therefore, with respect to the pseudo-orthonormal basis
		$
		\left\{\overline{e}_1,\overline{e}_2,\overline{e}_3,\overline{e}_4,\overline{e}_5\right\},
		$
		the matrix of $L$ is
		$$
		L_{\mu,s}
		=
		\left(
		\begin{array}{ccccc}
			\mu\cos s&-\mu\sin s&0&0&0\\
			\sin s&\cos s&0&0&0\\
			0&0&1&0&0\\
			0&0&0&1&0\\
			0&0&0&0&1
		\end{array}
		\right).
		$$
		If $\mu=1$, then
		$$
		L_{1,s}=R_s=
		\left(
		\begin{array}{ccccc}
			\cos s&-\sin s&0&0&0\\
			\sin s&\cos s&0&0&0\\
			0&0&1&0&0\\
			0&0&0&1&0\\
			0&0&0&0&1
		\end{array}
		\right).
		$$
		The transformation $R_s$ is the Euclidean rotation of angle $s$ in the spacelike plane
		$\Span\left\{\overline{e}_1,\overline{e}_2\right\}$, extended by the identity on
		$\Span\left\{\overline{e}_3,\overline{e}_4,\overline{e}_5\right\}$.
		
		If $\mu=-1$, then
		$$
		L_{-1,s}
		=
		\tilde{J}R_s
		=
		R_{-s}\tilde{J},
		$$
		where $\tilde{J}=\operatorname{diag}\left(-1,1,1,1,1\right)$.
		
		Thus, any isometry of $\mathbb{H}^4$ preserving $M$ is of the form $R_s$ or $\tilde{J}R_s$.
		
		Conversely, the matrices $R_s$ and $\tilde{J}R_s$ belong to $O^+\left(4,1\right)$. Moreover, we can check directly that
		$$
		R_s\Phi\left(u,t\right)
		=
		\Phi\left(u,t-s\right)
		$$
		and
		$$
		\left(\tilde{J}R_s\right)\Phi\left(u,t\right)
		=
		\Phi\left(u,-t+s\right).
		$$
		Hence, both $R_s$ and $\tilde{J}R_s$ preserve $M$. Therefore, these are all the isometries of $\mathbb{H}^4$ preserving the surface.
		
		Finally, we note that	
		$$
		R_s R_t=R_{s+t},
		\qquad
		R_0=I_5,
		$$
		Thus, the transformations $R_s$ form a one-parameter group. They act as Euclidean rotations in the spacelike plane $\Span\left\{\overline{e}_1,\overline{e}_2\right\}$ and fix pointwise its Lorentz-orthogonal complement $\Span\left\{\overline{e}_3,\overline{e}_4,\overline{e}_5\right\}$. For $s\notin2\pi\mathbb{Z}$, $R_s$ is nontrivial and fixes points in $\mathbb{H}^4$. More precisely, it fixes pointwise
		$\mathbb{H}^4\cap\Span\{\overline{e}_3,\overline{e}_4,\overline{e}_5\}$,
		which is a totally geodesic submanifold isometric to $\mathbb{H}^2$. Hence, the restriction of $R_s$ to $\mathbb{H}^4$ is an elliptic isometry.
		
		Furthermore, for any fixed $u$ and any $t_0\in\mathbb{R}$,
		$$
		\left\{
		R_s\Phi\left(u,t_0\right)
		\,\middle|\,
		s\in\mathbb{R}
		\right\}
		=
		\left\{
		\Phi\left(u,t\right)
		\,\middle|\,
		t\in\mathbb{R}
		\right\}.
		$$
		Thus, the curve $t\mapsto\Phi\left(u,t\right)$ is an orbit of the one-parameter group $\left\{R_s\, |\, s\in\mathbb{R}\right\}$.
	\end{proof}
	
	\subsubsection{The subcase $C<0$}
	%%%%%%%%%%%%%%%%%%%%%%%%%%%%%%%%%%%%%%%%%%%%%%%%%%%
	
	First, let us consider an integral curve $\gamma=\gamma\left(v\right)$ of $E_2$, viewed in $\mathbb{R}^5_1$. Then, as $\hat{\nabla}_{E_2}E_2=\hat{\kappa}\xi$ and $\hat{\nabla}_{E_2}\xi=\hat{\kappa}E_2$, it follows that, along $\gamma$, we have
	\begin{equation*}
		\left\{
		\begin{array}{l}
			\gamma''\left(v\right)=\hat{\kappa}\xi\left(v\right)\\
			\xi'\left(v\right)=\hat{\kappa}\gamma'\left(v\right)
		\end{array}
		\right..
	\end{equation*}
	Therefore $\gamma$ must satisfy
	$$
	\gamma'''\left(v\right)-\hat{\kappa}^2\gamma'\left(v\right)=0.
	$$
	A standard computation leads to
	$$
	\gamma\left(v\right)=\overline{A}_0+\frac{\sinh\left(\hat{\kappa}v\right)}{\hat{\kappa}}\overline{A}_1+\frac{\cosh\left(\hat{\kappa}v\right)}{\hat{\kappa}}\overline{A}_2,
	$$
	where $\overline{A}_0$, $\overline{A}_1$ and $\overline{A}_2$ are constant vectors in $\mathbb{R}_1^5$ 
	given by
	$$
	\overline{A}_0=\gamma\left(0\right)-\frac{1}{\hat{\kappa}}\xi\left(0\right),\qquad \overline{A}_1=\gamma'\left(0\right)=E_2\left(0\right), \qquad \overline{A}_2=\xi\left(0\right).
	$$
	Since $\langle\gamma', \gamma'\rangle=1$, it follows that
	$$
	\left\langle \overline{A}_1, \overline{A}_2\right\rangle=0 \qquad\text{and}\qquad \left\langle \overline{A}_1, \overline{A}_1\right\rangle=-\left\langle \overline{A}_2, \overline{A}_2\right\rangle=1.
	$$
	If we denote by 
	$$
	\overline{B}_0=\overline{A}_0, \qquad \overline{B}_1=\frac{1}{\hat{\kappa}}\overline{A}_1 \qquad \text{and}\qquad \overline{B}_2=\frac{1}{\hat{\kappa}} \overline{A}_2,  
	$$
	we rewrite
	\begin{equation*}
		\gamma\left(v\right)=\overline{B}_0+\sinh\left(\hat{\kappa} v\right)\overline{B}_1+\cosh\left(\hat{\kappa} v\right)\overline{B}_2,
	\end{equation*}
	where $\overline{B}_i \in \mathbb{R}_1^5$, $i=\overline{0,2}$ satisfy
	$$
	\overline{B}_0=\gamma\left(0\right)-\frac{1}{\hat{\kappa}}\xi\left(0\right),\qquad \overline{B}_1=\frac{1}{\hat{\kappa}}\gamma'\left(0\right)=\frac{1}{\hat{\kappa}}E_2\left(0\right), \qquad \overline{B}_2=\frac{1}{\hat{\kappa}}\xi\left(0\right),
	$$
	and
	$$
	\left\langle \overline{B}_1,\overline{B}_2\right\rangle=0 \qquad\text{and}\qquad \left\langle \overline{B}_1,\overline{B}_1\right\rangle=-\left\langle \overline{B}_2,\overline{B}_2\right\rangle=
	\frac{1}{\hat{\kappa}^2}.
	$$
	Thus, we can state the following properties of the integral curves of $E_2$.
	\begin{proposition}\label{integralcurveE2-case3}
		Let $\varphi:M^2 \to \mathbb{H}^4$ be a PNMC biconservative surface, and let $\gamma=\gamma\left(v\right)$ be an integral curve of $E_2$, viewed as a curve in $\mathbb{R}_1^5$. Assume that $\hat{\nabla}_{E_2}E_2$ is a non-zero timelike vector field. Then the following hold:
		\begin{itemize}
			\item[(i)] the curvature $\hat{\kappa}$ of $\gamma$ is a positive constant along $\gamma$;
			\item[(ii)] the vector field $\xi=\left(1/\hat{\kappa}\right) \hat{\nabla}_{E_2} E_2$ satisfies $\hat{\nabla}_{E_2}\xi=\hat{\kappa}E_2$;
			\item[(iii)] the curve $\gamma$ is given by
			$$
			\gamma\left(v\right)=\overline{B}_0+\sinh\left(\hat{\kappa}v\right)\overline{B}_1+\cosh\left(\hat{\kappa}v\right)\overline{B}_2,
			$$
			where
			$$
			\overline{B}_0=\gamma\left(0\right)-\frac{1}{\hat{\kappa}}\xi\left(0\right),\qquad
			\overline{B}_1=\frac{1}{\hat{\kappa}}E_2\left(0\right),\qquad
			\overline{B}_2=\frac{1}{\hat{\kappa}}\xi\left(0\right),
			$$
			and
			$$
			\left\langle \overline{B}_1,\overline{B}_2\right\rangle=0,\qquad
			\left\langle \overline{B}_1,\overline{B}_1\right\rangle=
			-\left\langle \overline{B}_2,\overline{B}_2\right\rangle=\frac{1}{\hat{\kappa}^2}.
			$$
			In particular, $\gamma$ is a branch of a hyperbola in the Lorentzian affine plane
			$$
			\overline{B}_0+\operatorname{span}\{\overline{B}_1,\overline{B}_2\}\subset\mathbb{R}_1^5.
			$$
		\end{itemize}
	\end{proposition}
	To obtain the parametrization of $M$ as a surface in $\mathbb{R}^5_1$, let $p_0 \in M$ be an arbitrarily fixed point of $M$ and $\sigma=\sigma\left(u\right)$ be an integral curve of $E_1$ with $\sigma\left(0\right)=p_0$. Consider $\left\{\phi_v\right\}_{v \in \mathbb{R}}$ the flow of $E_2$ near the point $p_0$. Then, for any $u \in\left(-\varepsilon, \varepsilon\right)$ and for any $v \in \mathbb{R}$, the parametrization $\Phi=\Phi\left(u, v\right)$ of $M$ in $\mathbb{R}^5_1$ is given by 
	$$
	\Phi\left(u, v\right)=\phi_{\sigma\left(u\right)}\left(v\right)=\overline{B}_0\left(u\right)+\sinh\left(\hat{\kappa}\left(u\right)v\right)\overline{B}_1\left(u\right)+\cosh \left(\hat{\kappa}\left(u\right)v\right)\overline{B}_2\left(u\right),
	$$
	where the vector-valued functions $\overline{B}_0, \overline{B}_1, \overline{B}_2$, which are uniquely determined by the surface, are given by
	\begin{equation*}
		\overline{B}_0\left(u\right)=\sigma\left(u\right)-\frac{1}{\hat{\kappa}\left(u\right)}\xi\left(u,0\right), \qquad \overline{B}_1\left(u\right)=\frac{1}{\hat{\kappa}\left(u\right)}E_2\left(u,0\right),\qquad \overline{B}_2\left(u\right)=\frac{1}{\hat{\kappa}\left(u\right)}\xi\left(u,0\right),
	\end{equation*}
	and they satisfy
	\begin{equation*}
		\left\langle \overline{B}_1\left(u\right), \overline{B}_2\left(u\right)\right\rangle=0\qquad\text{and}\qquad\left\langle \overline{B}_1\left(u\right), \overline{B}_1\left(u\right)\right\rangle=-\left\langle \overline{B}_2\left(u\right), \overline{B}_2\left(u\right)\right\rangle=\frac{1}{\hat{\kappa}^2\left(u\right)}.
	\end{equation*}
	Let us consider the vector-valued functions in $\mathbb{R}^5_1$ given by
	\begin{equation}\label{expressions-bi-C<0}
		\overline{b}_i\left(u\right)=\hat{\kappa}\left(u\right)\overline{B}_i\left(u\right),\qquad i=1,2.
	\end{equation}
	Then,
	$$
	\overline{b}_1\left(u\right)=E_2\left(u,0\right), \qquad \overline{b}_2\left(u\right)=\xi\left(u,0\right),
	$$
	and the parametrization can be expressed as
	$$
	\Phi\left(u, v\right)=\sigma\left(u\right)+\frac{1}{\hat{\kappa}\left(u\right)}\left(\sinh\left(\hat{\kappa}\left(u\right)v\right)\overline{b}_1\left(u\right)+\left(\cosh \left(\hat{\kappa}\left(u\right)v\right)-1\right)\overline{b}_2\left(u\right)\right).
	$$
	Clearly, $\Phi\left(u,0\right)=\sigma\left(u\right)$. Moreover, by Proposition \ref{Proposition 3.2} and the first relation in \eqref{hat-nabla-xi}, it follows that $E_2$ and $\xi$  are constant vector fields in $\mathbb{R}^5_1$ along $\sigma$. Therefore, the vector-valued functions $\overline{b}_i=\overline{b}_i\left(u\right)$ are constant in $\mathbb{R}^5_1$; hence, they can be identified with the constant vectors $\overline{b}_i\in\mathbb{R}^5_1$ satisfying
	$$
	\langle\overline{b}_1,\overline{b}_2\rangle=0\qquad \text{and}\qquad \langle\overline{b}_1,\overline{b}_1\rangle=-\langle\overline{b}_2,\overline{b}_2\rangle=1.
	$$
	In order to get a simpler expression of $\Phi$ we can consider the following change of coordinates
	$\left(u, v\right) \rightarrow\left(u, t=\hat{\kappa}\left(u\right) v\right)$. With respect to these new local coordinates, the parametrization $\Phi$ can be expressed as
	$$
	\Phi\left(u, t\right)=\sigma\left(u\right)+\frac{1}{\hat{\kappa}\left(u\right)}\left(\sinh t\ \overline{b}_1+\left(\cosh t-1\right)\overline{b}_2\right).
	$$
	From Proposition \ref{Proposition 3.2} we know that $E_2$ is constant along $\sigma$, so 
	$$
	\langle\overline{b}_1,\sigma\left(u\right)\rangle=0.
	$$ 
	Then, since $\langle \Phi\left(u,v\right),\Phi\left(u,v\right)\rangle=-1$, one obtains
	$$
	\langle \sigma\left(u\right),\overline{b}_2\rangle=-\frac{1}{\hat{\kappa}\left(u\right)}.
	$$
	Before stating the conclusion of this third part, we state the following result that can be obtained by straightforward computations.
	\begin{proposition}
		Let $\mathbb{H}^3$ be the totally geodesic hypersurface of $\mathbb{H}^4$ determined by the constant vector $\overline{b}_1\in\mathbb{R}^5_1$ given by \eqref{expressions-bi-C<0}. Then, the torsion of $\sigma$ satisfies \eqref{general-tau} in Lemma \ref{lemma-1}, and the constant vector $\overline{b}_2\in\mathbb{R}^5_1$ given by \eqref{expressions-bi-C<0} coincides with the vector $\overline{b}_2$ given by \eqref{expression-b2} in Lemma \ref{lemma-2}, for
		$$
		\delta=-1 \qquad\text{and}\qquad \Theta\left(u\right)=-\frac{1}{3\sqrt{-C}f^{3/4}\left(u\right)}.
		$$
	\end{proposition}
	We are now ready to state the main result of this third case.
	\begin{theorem}\label{theorem-C<0}
		Let $\varphi:M^2\rightarrow \mathbb{H}^4$ be a PNMC biconservative immersion and denote $\Phi=\iota \circ \varphi: M \rightarrow \mathbb{R}_1^5$, where $\iota: \mathbb{H}^4 \rightarrow \mathbb{R}_1^5$ is the canonical inclusion. Identifying $M$ with its image, the surface $M$ can be locally parametrized as
		$$
		\Phi\left(u, t\right)=\sigma\left(u\right)+\frac{1}{\hat{\kappa}\left(u\right)}\left(\sinh t\ \overline{b}_1+\left(\cosh t-1\right)\overline{b}_2\right),
		$$
		where
		\begin{itemize}
			\item[(i)] $f=f\left(u\right)$ is a positive solution of the fist-order ordinary differential equation 
			$$
			f'=\frac{4}{3}f\sqrt{1+9Cf^{3/2}-9f^2-c^2f^3}>0,
			$$
			where $c$ is a real constant and $C$ is a negative real constant;
			\item[(ii)]$\hat{\kappa}=\hat{\kappa}\left(u\right)$ is the positive function given by
			$$
			\hat{\kappa}\left(u\right)=3\sqrt{-C}f^{3/4}\left(u\right);
			$$
			\item[(iii)] $\overline{b}_1$ is a constant vector in $\mathbb{R}_1^5$ such that 
			$\langle\overline{b}_1,\overline{b}_1\rangle=1$;
			\item[(iv)] $\sigma=\sigma\left(u\right)$ is a curve parametrized by arc-length in $\mathbb{R}_1^5$, which lies in a totally geodesic hypersurface $\mathbb{H}^3=\mathbb{H}^4 \cap \Pi$; the hyperplane $\Pi$ contains the origin and is orthogonal to $\overline{b}_1$. Moreover, the curvature and torsion of $\sigma$, as a curve in $\mathbb{H}^3$, are given in Proposition \ref{Proposition 3.2} and $\overline{b}_2$ is a constant vector in $\mathbb{R}_1^5$ which satisfies
			$$
			\langle\overline{b}_2,\overline{b}_1\rangle=0, \qquad \langle\overline{b}_2,\overline{b}_2\rangle=-1 \qquad \text{and}\qquad  \langle\sigma\left(u\right), \overline{b}_2\rangle=-\frac{1}{\hat{\kappa}\left(u\right)}.
			$$
		\end{itemize}
	\end{theorem}
	Further, we give some geometric properties of the parametrization given in Theorem \ref{theorem-C<0}.
	\begin{proposition}
		Let $\left\{\overline{e}_1=\overline{b}_1,\overline{e}_2,\overline{e}_3,\overline{e}_4,\overline{e}_5=\overline{b}_2\right\}$ be a pseudo-orthonormal basis in $\mathbb{R}_1^5$. Then, the isometries of $\mathbb{H}^4$ preserving the surface $M$ are precisely the transformations $B_s$ and $\tilde{J}B_s$, with $s\in\mathbb{R}$, where
		$$
		B_s=
		\left(
		\begin{array}{ccccc}
			\cosh s&0&0&0&\sinh s\\
			0&1&0&0&0\\
			0&0&1&0&0\\
			0&0&0&1&0\\
			\sinh s&0&0&0&\cosh s
		\end{array}
		\right)
		$$
		and $\tilde{J}=\operatorname{diag}\left(-1,1,1,1,1\right)$. The transformations $B_s$ form a one-parameter group of Lorentz boosts in the Lorentzian plane $\Span\left\{\overline{e}_1,\overline{e}_5\right\}$. More precisely,
		$$
		B_s\Phi\left(u,t\right)
		=
		\Phi\left(u,t+s\right),
		$$
		and
		$$
		\left(\tilde{J}B_s\right)\Phi\left(u,t\right)
		=
		\Phi\left(u,-t-s\right).
		$$
		In particular, for any fixed $u$, the curve $t\mapsto\Phi\left(u,t\right)$ is an orbit of the one-parameter group $\left\{B_s\, |\, s\in\mathbb{R}\right\}$.
	\end{proposition}
	\begin{proof}
		With respect to the above pseudo-orthonormal basis, the parametrization becomes
		$$
		\Phi\left(u,t\right)
		=
		\sigma\left(u\right)
		+
		\frac{1}{\hat{\kappa}\left(u\right)}
		\left(
		\sinh t\,\overline{e}_1
		+
		\left(\cosh t-1\right)\overline{e}_5
		\right).
		$$
		Furthermore, we have
		$$
		\langle\sigma,\overline{e}_1\rangle=0,
		\qquad
		\langle\sigma,\overline{e}_5\rangle
		=
		-\frac{1}{\hat{\kappa}\left(u\right)}.
		$$
		If we denote
		\begin{equation}\label{alpha-sigma-C<0}
			\alpha\left(u\right)
			=
			\sigma\left(u\right)
			-
			\frac{1}{\hat{\kappa}\left(u\right)}\overline{e}_5,
		\end{equation}
		then
		$$
		\alpha\left(u\right)
		\in
		\Span\{\overline{e}_2,\overline{e}_3,\overline{e}_4\},
		$$
		and the parametrization can be written as
		$$
		\Phi\left(u,t\right)
		=
		\alpha\left(u\right)
		+
		\frac{1}{\hat{\kappa}\left(u\right)}
		\left(
		\sinh t\,\overline{e}_1
		+
		\cosh t\,\overline{e}_5
		\right).
		$$
		Thus, for any fixed $u$, the curve $t\mapsto\Phi\left(u,t\right)$ is one branch of a hyperbola lying in the affine Lorentzian plane
		$\alpha\left(u\right)+\Span\{\overline{e}_1,\overline{e}_5\}$.
		
		Let $L\in O^+\left(4,1\right)$ be an isometry of $\mathbb{H}^4$ preserving $M$. Using the same argument as in the proof of Proposition \ref{prop-C=0}, $L$ preserves each curve $t\mapsto\Phi\left(u,t\right)$. Hence,
		\begin{equation}\label{L-phi-dif-C<0}
			L\Phi\left(u,t\right)=
			\Phi\left(u,\psi_u\left(t\right)\right),
		\end{equation}
		where $\psi_u$ represents a reparametrization of the hyperbola $t\mapsto\Phi\left(u,t\right)$, for an arbitrarily fixed $u$.
		Moreover, we can see that
		$$
		\Phi_t
		=
		\frac{1}{\hat{\kappa}\left(u\right)}
		\left(
		\cosh t\,\overline{e}_1
		+
		\sinh t\,\overline{e}_5
		\right),
		$$
		and
		$$
		\langle\Phi_t,\Phi_t\rangle
		=
		\frac{1}{\hat{\kappa}^2\left(u\right)}.
		$$
		Differentiating \eqref{L-phi-dif-C<0} with respect to $t$ and using that $L$ is an isometry of $\mathbb{H}^4$, we obtain
		$$
		\frac{1}{\hat{\kappa}^2\left(u\right)}
		=
		\left\langle
		L\Phi_t\left(u,t\right),
		L\Phi_t\left(u,t\right)
		\right\rangle
		=
		\left(\psi_u'\left(t\right)\right)^2
		\left\langle
		\Phi_t\left(u,\psi_u\left(t\right)\right),
		\Phi_t\left(u,\psi_u\left(t\right)\right)
		\right\rangle
		=
		\left(\psi_u'\left(t\right)\right)^2
		\frac{1}{\hat{\kappa}^2\left(u\right)}.
		$$
		Therefore,
		$$
		\psi_u\left(t\right)
		=
		\mu\left(u\right)t+a\left(u\right),
		$$
		where
		$$
		\mu\left(u\right)=\mu\in\{-1,1\},
		\qquad
		a\left(u\right)\in\mathbb{R}.
		$$
		Hence,
		\begin{equation}\label{L-phi(u,t)-C<0}
			L\Phi\left(u,t\right)
			=
			\Phi\left(u,\mu t+a\left(u\right)\right),
			\qquad
			\mu\in\{-1,1\}.
		\end{equation}
		Next, we will determine $a\left(u\right)$ and $L$. Expanding both sides of \eqref{L-phi(u,t)-C<0}, we obtain
		\begin{equation}\label{extended-L-phi-C<0}
			L\alpha\left(u\right)
			+
			\frac{1}{\hat{\kappa}\left(u\right)}
			\left(
			\sinh t\,L\overline{e}_1
			+
			\cosh t\,L\overline{e}_5
			\right)
			=
			\alpha\left(u\right)
			+
			\frac{1}{\hat{\kappa}\left(u\right)}
			\left(
			\sinh\left(\mu t+a\left(u\right)\right)\,\overline{e}_1
			+
			\cosh\left(\mu t+a\left(u\right)\right)\,\overline{e}_5
			\right).
		\end{equation}
		Since $\mu\in\{-1,1\}$, we have
		$$
		\sinh\left(\mu t+a\left(u\right)\right)
		=
		\mu\cosh a\left(u\right)\sinh t
		+
		\sinh a\left(u\right)\cosh t
		$$
		and
		$$
		\cosh\left(\mu t+a\left(u\right)\right)
		=
		\mu\sinh a\left(u\right)\sinh t
		+
		\cosh a\left(u\right)\cosh t.
		$$
		Therefore, comparing the coefficients of $\sinh t$, $\cosh t$, and the constant terms in \eqref{extended-L-phi-C<0}, we obtain
		$$
		L\overline{e}_1
		=
		\mu\cosh a\left(u\right)\,\overline{e}_1
		+
		\mu\sinh a\left(u\right)\,\overline{e}_5,
		$$
		$$
		L\overline{e}_5
		=
		\sinh a\left(u\right)\,\overline{e}_1
		+
		\cosh a\left(u\right)\,\overline{e}_5,
		$$
		and
		$$
		L\alpha\left(u\right)=\alpha\left(u\right).
		$$
		Since $\overline{e}_i$ and $L\overline{e}_i$ are independent of $u$, both $\sinh a\left(u\right)$ and $\cosh a\left(u\right)$ are independent of $u$. Hence, $a\left(u\right)$ is independent of $u$. Thus, there exists a real constant $s$ such that
		$$
		a\left(u\right)=\mu s.
		$$
		Consequently,
		$$
		L\Phi\left(u,t\right)
		=
		\Phi\left(u,\mu\left(t+s\right)\right),
		$$
		\begin{equation}\label{Le1-C<0}
			L\overline{e}_1
			=
			\mu\cosh s\,\overline{e}_1
			+
			\sinh s\,\overline{e}_5,
		\end{equation}
		and
		\begin{equation}\label{Le5-C<0}
			L\overline{e}_5
			=
			\mu\sinh s\,\overline{e}_1
			+
			\cosh s\,\overline{e}_5.
		\end{equation}
		Hence,
		$$
		L\left(
		\Span\{\overline{e}_1,\overline{e}_5\}
		\right)
		=
		\Span\{\overline{e}_1,\overline{e}_5\}.
		$$
		Since $L$ preserves the Lorentzian scalar product, it also preserves the Lorentz-orthogonal complement of the plane
		$\Span\left\{\overline{e}_1,\overline{e}_5\right\}$. Therefore,
		$$
		L\left(
		\Span\{\overline{e}_2,\overline{e}_3,\overline{e}_4\}
		\right)
		=
		\Span\{\overline{e}_2,\overline{e}_3,\overline{e}_4\}.
		$$
		Further, we recall that, as $\sigma$ lies full in $\mathbb{H}^3$, which is determined by $\overline{e}_1^\perp$, we have
		$$
		\Span\{\sigma\left(u\right)\, |\, u\in I\}
		=
		\overline{e}_1^\perp
		=
		\Span\{\overline{e}_2,\overline{e}_3,\overline{e}_4,\overline{e}_5\}.
		$$
		Then, as $\alpha\left(u\right)$ is given in \eqref{alpha-sigma-C<0} and the $\overline{e}_5$-component of $\alpha\left(u\right)$ vanishes we obtain that 
		$$
		\Span 
		\left\{\alpha\left(u\right)\, |\,  u\in I \right\}
		=
		\Span
		\{\overline{e}_2,\overline{e}_3,\overline{e}_4\}.
		$$
		Now, since $L$ fixes any vector $\alpha\left(u\right)$, we obtain
		$$
		L\overline{e}_2=\overline{e}_2,
		\qquad
		L\overline{e}_3=\overline{e}_3,
		\qquad
		L\overline{e}_4=\overline{e}_4.
		$$
		Therefore, with respect to the pseudo-orthonormal basis
		$\left\{\overline{e}_1,\overline{e}_2,\overline{e}_3,\overline{e}_4,\overline{e}_5\right\}$,
		the matrix of $L$ is
		$$
		L_{\mu,s}
		=
		\left(
		\begin{array}{ccccc}
			\mu\cosh s&0&0&0&\mu\sinh s\\
			0&1&0&0&0\\
			0&0&1&0&0\\
			0&0&0&1&0\\
			\sinh s&0&0&0&\cosh s
		\end{array}
		\right).
		$$
		If $\mu=1$, then
		$$
		L_{1,s}=B_s=
		\left(
		\begin{array}{ccccc}
			\cosh s&0&0&0&\sinh s\\
			0&1&0&0&0\\
			0&0&1&0&0\\
			0&0&0&1&0\\
			\sinh s&0&0&0&\cosh s
		\end{array}
		\right).
		$$
		The transformation $B_s$ is the Lorentz boost of rapidity $s$ in the Lorentzian plane
		$\Span\left\{\overline{e}_1,\overline{e}_5\right\}$, extended by the identity on
		$\Span\left\{\overline{e}_2,\overline{e}_3,\overline{e}_4\right\}$.
		
		If $\mu=-1$, then
		$$
		L_{-1,s}
		=
		\tilde{J}B_s
		=
		B_{-s}\tilde{J},
		$$
		where
		$\tilde{J}=\operatorname{diag}\left(-1,1,1,1,1\right)$. 
		
		Thus, any isometry of $\mathbb{H}^4$ preserving $M$ is of the form $B_s$ or $\tilde{J}B_s$.
		
		Conversely, the matrices $B_s$ and $\tilde{J}B_s$ belong to $O^+\left(4,1\right)$. Moreover, we can check directly that
		$$
		B_s\Phi\left(u,t\right)
		=
		\Phi\left(u,t+s\right)
		$$
		and
		$$
		\left(\tilde{J}B_s\right)\Phi\left(u,t\right)
		=
		\Phi\left(u,-t-s\right).
		$$
		Hence, both $B_s$ and $\tilde{J}B_s$ preserve $M$. Therefore, these are all the isometries of $\mathbb{H}^4$ preserving the surface.
		
		Finally, we note that
		$$
		B_sB_t=B_{s+t},
		\qquad
		B_0=I_5.
		$$
		Thus, the transformations $B_s$ form a one-parameter group. They act as Lorentz boosts in the Lorentzian plane
		$\Span\left\{\overline{e}_1,\overline{e}_5\right\}$ and fix pointwise its Lorentz-orthogonal complement
		$\Span\left\{\overline{e}_2,\overline{e}_3,\overline{e}_4\right\}$.
		
		For $s\neq0$, we first observe that $B_s$ has no fixed point in $\mathbb{H}^4$. Indeed, if $B_s\overline{x}=\overline{x}$,
		then $\overline{x}$ belongs to the eigenspace of $B_s$ corresponding to the eigenvalue $1$. Since $B_s$ acts as the identity on
		$\Span\{\overline{e}_2,\overline{e}_3,\overline{e}_4\}$
		and its restriction to $\Span\{\overline{e}_1,\overline{e}_5\}$ has no eigenvalue equal to $1$ for $s\neq0$, we have
		$$
		\ker\left(B_s-I_5\right)
		=
		\Span\{\overline{e}_2,\overline{e}_3,\overline{e}_4\}.
		$$
		Any nonzero vector $\overline{x}$ in this subspace satisfies
		$\langle\overline{x},\overline{x}\rangle>0$. Therefore, $\overline{x}\notin\mathbb{H}^4$, and hence $B_s$ has no fixed point in $\mathbb{H}^4$.
		
		We now determine the invariant null directions of $B_s$. A one-dimensional subspace $\mathbb{R}\overline{x}$ is invariant under $B_s$ if and only if $B_s\overline{x}=\lambda\overline{x}$, for some $\lambda\in\mathbb{R}$. Thus, the invariant directions of $B_s$ are precisely its eigendirections.
		
		Since $B_s$ acts as the identity on $\Span\{\overline{e}_2,\overline{e}_3,\overline{e}_4\}$,
		and its restriction to the Lorentzian plane $\Span\{\overline{e}_1,\overline{e}_5\}$ is represented by
		$$
		\left(
		\begin{array}{cc}
			\cosh s&\sinh s\\
			\sinh s&\cosh s
		\end{array}
		\right),
		$$
		we can check that the eigenvalues of $B_s$ are
		$1$, $e^s$ and $e^{-s}$, with multiplicities $3$, $1$, and $1$, respectively. The corresponding eigenspaces are
		$$
		E_1
		=
		\Span\{\overline{e}_2,\overline{e}_3,\overline{e}_4\},
		$$
		$$
		E_{e^s}
		=
		\Span\{\overline{e}_1+\overline{e}_5\},
		$$
		and
		$$
		E_{e^{-s}}
		=
		\Span\{-\overline{e}_1+\overline{e}_5\}.
		$$
		Indeed,
		$$
		B_s\left(\overline{e}_1+\overline{e}_5\right)
		=
		e^s\left(\overline{e}_1+\overline{e}_5\right)
		$$
		and
		$$
		B_s\left(-\overline{e}_1+\overline{e}_5\right)
		=
		e^{-s}\left(-\overline{e}_1+\overline{e}_5\right).
		$$
		Moreover,
		$$
		\langle
		\overline{e}_1+\overline{e}_5,
		\overline{e}_1+\overline{e}_5
		\rangle
		=
		0
		$$
		and
		$$
		\langle
		-\overline{e}_1+\overline{e}_5,
		-\overline{e}_1+\overline{e}_5
		\rangle
		=
		0,
		$$
		so these two eigendirections are null.
		
		On the other hand, the eigenspace $E_1$ is positive definite and therefore contains no nonzero null vector. Since every invariant direction of $B_s$ is an eigendirection, it follows that
		$$
		\mathbb{R}\left(\overline{e}_1+\overline{e}_5\right)
		\qquad\text{and}\qquad
		\mathbb{R}\left(-\overline{e}_1+\overline{e}_5\right)
		$$
		are the only invariant null directions of $B_s$.
		
		Therefore, the restriction of $B_s$ to $\mathbb{H}^4$ has no fixed point in $\mathbb{H}^4$ and has exactly two fixed points at infinity, corresponding to these two null directions. Hence, for $s\neq0$, its restriction to $\mathbb{H}^4$ is a hyperbolic isometry.
		
		Furthermore, for any fixed $u$ and any $t_0\in\mathbb{R}$,
		$$
		\left\{
		B_s\Phi\left(u,t_0\right)
		\,\middle|\,
		s\in\mathbb{R}
		\right\}
		=
		\left\{
		\Phi\left(u,t\right)
		\,\middle|\,
		t\in\mathbb{R}
		\right\}.
		$$
		Thus, the curve $t\mapsto\Phi\left(u,t\right)$ is an orbit of the one-parameter group $\left\{B_s\, |\, s\in\mathbb{R}\right\}$.
	\end{proof}
We end this paper with some geometric properties of the directrix $\sigma$. These properties are given in terms of the causal character of the constant vector $\overline{b}_2$ associated with $\sigma\subset\mathbb{H}^3$.
	
If $C=0$, then $\overline{b}_2$ is null and
		$$
		\langle\sigma\left(u\right),\overline{b}_2\rangle
		=
		-\frac{1}{2f^{3/4}\left(u\right)}.
		$$
		The null direction $\mathbb{R}\overline{b}_2$ determines a point at infinity of $\mathbb{H}^3$, denoted by $[\overline{b}_2]$, and, in $\overline{b}_1^\perp$, the negative level sets of the function $\overline{x}\mapsto\langle \overline{x},\overline{b}_2\rangle$ intersected with $\mathbb{H}^3$ are horospheres centered at this point at infinity. Since
		$$
		\langle \sigma'\left(u\right),\overline{b}_2\rangle
		=
		\frac{3}{8}f^{-7/4}\left(u\right)f'\left(u\right)>0,
		$$
		the curve $\sigma$ is nowhere tangent to the corresponding negative level sets and therefore passes transversely through the one-parameter family of horospheres centered at the same point at infinity $[\overline{b}_2]$. Moreover, since $f$ is strictly increasing, the negative level parameter
		$$
		-\frac{1}{2f^{3/4}\left(u\right)}
		$$
		is strictly increasing. Hence, as $u$ increases, the corresponding horospheres move toward their common center $[\overline{b}_2]$.
		
		If $C>0$, then $\overline{b}_2$ is spacelike and
		$$
		\langle\sigma\left(u\right),\overline{b}_2\rangle
		=
		\frac{1}{\hat{\kappa}\left(u\right)}.
		$$
		We note that $\mathbb{H}^2=\mathbb{H}^3\cap\overline{b}_2^\perp$ is a totally geodesic surface in $\mathbb{H}^3$; here, $\overline{b}_2^\perp$ is the orthogonal complement of $\overline{b}_2$ in $\overline{b}_1^\perp$. Also, we can see that the positive level sets
		$$
		\left\{
		\overline{x}\in\mathbb{H}^3\, |\,
		\langle \overline{x},\overline{b}_2\rangle=\lambda
		\right\}
		$$	
		are equidistant surfaces with respect to $\mathbb{H}^2$. More precisely, as $\lambda>0$, the distance from such a level set to $\mathbb{H}^2$ is $\operatorname{arsinh}\left(\lambda\right)$. Indeed, for an arbitrary fixed $\lambda>0$, let us consider a point $\overline{x}\in\mathbb{H}^3$ such that $\langle\overline{x},\overline{b}_2\rangle=\lambda$. Then, $\overline{x}$ can be written as
		$\overline{x}=\overline{v}+\lambda\overline{b}_2$, where $\overline{v}\in\overline{b}_2^\perp$. Using $\langle\overline{x},\overline{x}\rangle=-1$, we obtain
		$$
		\langle\overline{v},\overline{v}\rangle
		=
		-\left(1+\lambda^2\right).
		$$
		Further, by standard computations, as $\langle\overline{x},\overline{v}\rangle<0$, we can prove that
		$$
		\overline{y}
		=
		\frac{\overline{v}}{\sqrt{1+\lambda^2}}
		$$
		belongs to the upper part of the hyperboloid, i.e., $\mathbb{H}^2$, and therefore
		$$
		\overline{x}
		=
		\sqrt{1+\lambda^2}\,\overline{y}
		+
		\lambda\overline{b}_2.
		$$
		Writing $\lambda=\sinh\rho$, we obtain
		$$
		\overline{x}
		=
		\cosh\rho\,\overline{y}
		+
		\sinh\rho\,\overline{b}_2.
		$$
		We note that the curve
		$$
		\gamma\left(t\right)
		=
		\cosh t\,\overline{y}
		+
		\sinh t\,\overline{b}_2
		$$
		is an arc-length parametrized geodesic of $\mathbb{H}^3$ with $\gamma\left(0\right)=\overline{y}\in\mathbb{H}^2$, $\gamma\left(\rho\right)=\overline{x}$ and $\gamma'\left(0\right)=\overline{b}_2$. Consequently,
		$$
		d\left(\overline{x},\mathbb{H}^2\right)
		=
		\rho
		=
		\operatorname{arsinh}\left(\lambda\right).
		$$
		As the point $\overline{x}$ was arbitrarily fixed, the positive level sets $\left\{
		\overline{x}\in\mathbb{H}^3\, |\,
		\langle \overline{x},\overline{b}_2\rangle=\lambda
		\right\}$ are equidistant surfaces with respect to the totally geodesic surface $\mathbb{H}^2$.
		
		In particular, since
		$$
		\langle\sigma\left(u\right),\overline{b}_2\rangle
		=
		\frac{1}{\hat{\kappa}\left(u\right)},
		$$
		we have
		$$
		d\left(\sigma\left(u\right),\mathbb{H}^2\right)
		=
		\operatorname{arsinh}
		\left(
		\frac{1}{\hat{\kappa}\left(u\right)}
		\right).
		$$
		Moreover, since the level parameter $1/\hat{\kappa}\left(u\right)$ is strictly decreasing, the curve $\sigma$ is nowhere tangent to the corresponding level sets and therefore passes transversely through this one-parameter family of equidistant surfaces. Also, we can see that when $u$ increases the corresponding equidistant surfaces move towards the totally geodesic surface $\mathbb{H}^2$.
		
		If $C<0$, then $\overline{b}_2$ is timelike and, using Remark \ref{remark-b2-H3}, it follows that $\overline{b}_2\in\mathbb{H}^3$. Moreover,
		$$
		\langle\sigma\left(u\right),\overline{b}_2\rangle
		=
		-\frac{1}{\hat{\kappa}\left(u\right)}.
		$$
		Since
		$$
		f'
		=
		\frac{4}{3}f
		\sqrt{1+9Cf^{3/2}-9f^2-c^2f^3}>0,
		$$
		it follows that
		$$
		1+9Cf^{3/2}-9f^2-c^2f^3>0.
		$$
		Since $C<0$, we obtain $9\left(-C\right)f^{3/2}<1$, and therefore, as $
		\hat{\kappa}\left(u\right)=3\sqrt{-C}\,f^{3/4}\left(u\right)$, it follows that
		$0<\hat{\kappa}\left(u\right)<1$. Hence,
		$$
		\frac{1}{\hat{\kappa}\left(u\right)}>1.
		$$
		Now, let us consider the level sets
		$$
		\left\{
		\overline{x}\in\mathbb{H}^3\,\big|\,
		\langle\overline{x},\overline{b}_2\rangle=-\lambda, \qquad \lambda>1
		\right\}.
		$$
		We can see that they are geodesic spheres centered at $\overline{b}_2$. More precisely, their radius is $\operatorname{arcosh}\left(\lambda\right)$. Indeed, for an arbitrary fixed $\lambda>1$, let $\overline{x}\in\mathbb{H}^3$ such that
		$\langle\overline{x},\overline{b}_2\rangle=-\lambda$.
		Since $\langle\overline{b}_2,\overline{b}_2\rangle=-1$, the vector $\overline{x}$ can be written as
		$\overline{x}=\overline{v}+\lambda\overline{b}_2$,
		where $\overline{v}\in\overline{b}_2^\perp$.
		Using $\langle\overline{x},\overline{x}\rangle=-1$, we obtain
		$$
		\langle\overline{v},\overline{v}\rangle
		=
		\lambda^2-1.
		$$
		Further, by a standard computation, we can prove that
		$$
		\overline{y}
		=
		\frac{\overline{v}}{\sqrt{\lambda^2-1}}
		$$
		is a unit spacelike vector orthogonal to $\overline{b}_2$.
		Consequently,
		$$
		\overline{x}
		=
		\lambda\overline{b}_2
		+
		\sqrt{\lambda^2-1}\,\overline{y}.
		$$
		Writing $\lambda=\cosh\rho$, with $\rho>0$, we obtain
		$$
		\overline{x}
		=
		\cosh\rho\,\overline{b}_2
		+
		\sinh\rho\,\overline{y}.
		$$
		We note that the curve
		$$
		\gamma\left(t\right)
		=
		\cosh t\,\overline{b}_2
		+
		\sinh t\,\overline{y}
		$$
		is an arc-length parametrized geodesic of $\mathbb{H}^3$ with
		$\gamma\left(0\right)=\overline{b}_2$, $\gamma\left(\rho\right)=\overline{x}$ and $\gamma'\left(0\right)=\overline{y}$. Consequently,
		$$
		d\left(\overline{x},\overline{b}_2\right)
		=
		\rho
		=
		\operatorname{arcosh}\left(\lambda\right).
		$$
		As the point $\overline{x}$ was arbitrarily fixed, the level set
		$$
		\left\{
		\overline{x}\in\mathbb{H}^3\,\big|\,
		\langle\overline{x},\overline{b}_2\rangle=-\lambda
		\right\}
		$$
		is the geodesic sphere centered at $\overline{b}_2$ with radius
		$\operatorname{arcosh}\left(\lambda\right)$.
		
		In particular, since
		$$
		\langle\sigma\left(u\right),\overline{b}_2\rangle
		=
		-\frac{1}{\hat{\kappa}\left(u\right)},
		$$
		we have
		$$
		d\left(\sigma\left(u\right),\overline{b}_2\right)
		=
		\operatorname{arcosh}
		\left(
		\frac{1}{\hat{\kappa}\left(u\right)}
		\right).
		$$
		Moreover, since the level parameter $1/\hat{\kappa}\left(u\right)$ is strictly decreasing, the curve $\sigma$ is nowhere tangent to the corresponding level sets and therefore passes transversely through this one-parameter family of geodesic spheres having the same center $\overline{b}_2$. Also, as $u$ increases, the radii of the corresponding geodesic spheres decrease, and hence these spheres move towards their common center $\overline{b}_2$.
		
		Thus, the three possibilities $C=0$, $C>0$, and $C<0$ correspond, respectively, to the null, spacelike, and timelike character of $\overline{b}_2$, and hence to horospherical, equidistant, and spherical geometries in $\mathbb{H}^3$.
	\begin{remark}
		One can check directly that each of the local parametrizations given in Theorem \ref{theorem-C=0}, Theorem \ref{theorem-C>0} and Theorem \ref{theorem-C<0}, defines a PNMC biconservative immersion from $M^2$ in $\mathbb{H}^4$, corresponding to the cases $C=0$, $C>0$ and $C<0$, respectively. In particular, in each case, the mean curvature function is indeed $f=f\left(u\right)$ and, up to isometries of the domain, the metric induced by $\Phi$ is given by
		$$
		g\left(u,v\right)=du^2+\frac{1}{f^{3/2}\left(u\right)}dv^2,
		$$
		which agrees with $\left(3.7\right)$ in \cite{SNMR2024}.
	\end{remark}


\begin{thebibliography}{0} 
		
		\bibitem{AN} \c{S}. Andronic, S. Nistor, \textit{Gap results for biharmonic submanifolds in spheres}, J. Math. Anal. Appl. 548 (2025), no. 1, Paper No. 129378, 34 pp.
		
		\bibitem{AK} \c{S}. Andronic,  A. Kayhan, \textit{Rigidity results for compact biconservative hypersurfaces in space forms}, J. Geom. Phys. 212 (2025), Paper No. 105460, 15 pp.
		
		\bibitem{AMOS} \c{S}. Andronic,  S. Montaldo, C. Oniciuc, A. Sanna, \textit{Biconservative Weingarten surfaces with flat normal bundle in  $N^4\left(\varepsilon\right)$}, Math. Nachr. 299 (2026), no. 7, 1624--1664.
		
		\bibitem{CMOP} R. Caddeo, S. Montaldo, C. Oniciuc and P. Piu, \textit{ Surfaces in the three-dimensional space forms with divergence-free stress-bienergy tensor}, Ann. Mat. Pura Appl. 193 (2014), 529--550.	
		
		\bibitem{F2025} D. Fetcu, \textit{The rigidity of biconservative surfaces in $\Sol_3$}, Ann. Mat. Pura Appl. (4) 204 (2025), no. 6, 2363--2375.
		
		\bibitem{FLO2021} D. Fetcu, E. Loubeau and C. Oniciuc, \textit{Bochner-Simons formulas and the rigidity of biharmonic submanifolds} J. Geom. Anal. 31 (2021), no. 2, 1732--1755.
		
		\bibitem{FO2022} D. Fetcu and C. Oniciuc, \textit{Biharmonic and biconservative hypersurfaces in space forms}, Differential geometry and global analysis—in honor of Tadashi Nagano, 65--90, Contemp. Math., 777, Amer. Math. Soc., 2022.
		
		\bibitem{FOP2015} D. Fetcu, C. Oniciuc and A. L. Pinheiro, \textit{CMC biconservative surfaces in $\mathbb{S}^n\times \mathbb{R}$ and $\mathbb{H}^n\times\mathbb{R}$}, J. Math. Anal. Appl. 425 (2015), no.1, 588--609.
		
		\bibitem{F2015} Y. Fu, \textit{Explicit classification of biconservative surfaces in Lorenz $3$-space forms}, Ann. Mat. Pura Appl. (4) 194 (2015), no.3, 805--822.
		
		\bibitem{Kayhan2025} A. Kayhan, \textit{Complete spacelike biconservative hypersurfaces in the Sitter space}, arXiv:2506.04744, 2025.
		
		\bibitem{Hasanis-Vlachos} T. Hasanis and T. Vlachos, \textit{Hypersurfaces in $\mathbb{E}^4$ with harmonic mean curvature vector field}, Math. Nachr. 172 (1995), 145--169.
		
		\bibitem{Jiang87} G.Y. Jiang, \textit{The conservation law for $2$-harmonic maps between Riemannian manifolds}, Acta Math. Sinica 30 (1987),  220--225.
		
		\bibitem{LMO} E. Loubeau, S. Montaldo and C. Oniciuc, \textit{The stress-energy tensor for biharmonic maps}, Math. Z. 259 (2008), 503--524.
		
		\bibitem{MOP2023} S. Montaldo, C. Oniciuc and A. Pampano, \textit{Closed biconservative hypersurfaces in spheres}, J. Math. Anal. Appl. 518 (2023), no.1, Paper No. 126697, 16 pp.
		
		\bibitem{MOR2016JGA} S. Montaldo, C. Oniciuc and A. Ratto, \textit{Biconservative surfaces}, J. Geom. Anal. 26 (2016), 313--329.
		
		\bibitem{NPhD17} S. Nistor, \textit{Biharmonicity and biconservativity topics in the theory of submanifolds}, PhD Thesis, 2017, doi: 10.13140/RG.2.2.27179.05924.
		
		\bibitem{NOTYS} S. Nistor, C. Oniciuc, N.C. Turgay and R. Ye\u gin \c Sen, \textit{Biconservative surfaces in the $4$-dimensional Euclidean sphere}, Ann. Mat. Pura Appl. 202 (2023), no. 5, 2345--2377. 
		
		\bibitem{SNMR2024} S. Nistor, M.Rusu,  \textit{Intrinsic characterization of biconservative surfaces in the $4$-dimensional hyperbolic space}, Mediterr. J. Math., 21 (2024), Paper no. 225, 27 pp.
		
		\bibitem{OC-20} Y.-L. Ou, B.-Y. Chen, \textit{Biharmonic Submanifolds and Biharmonic Maps in Riemannian Geometry}, World Scientific Publishing, Hackensack, NJ, 2020.
		
		\bibitem{Turgay2015} N.C. Turgay, \textit{H-hypersurfaces with three distinct principal curvatures in the Euclidean spaces}, Ann. Mat. Pura Appl. (4) 194 (2015), no. 6, 1795--1807.
		
		\bibitem{TurgayYegin2025} N.C. Turgay and R. Ye\u gin \c Sen, \textit{Biconservative surfaces in Robertson-Walker spaces}, Results Math. 80 (2025), no.3, Paper No. 77, 25 pp.
		
		
		\bibitem{YeginTurgay2018} N.C. Turgay and R. Ye\u gin \c Sen, \textit{On biconservative surfaces in $4$-dimensional Euclidean space}, J. Math. Anal. Appl. 460 (2018), 565--581.
		
		\bibitem{YeginTurgay2025} N.C. Turgay and R. Ye\u gin \c Sen, \textit{Biconservative PNMCV surfaces in the arbitrary dimensional Minkowski space}, J. Korean Math. Soc. 62 (2025), no. 1, 145--163.
		
	\end{thebibliography}
\end{document}